\documentclass[letterpaper,12pt]{article}
\usepackage{amsmath, amsthm}
\usepackage{amssymb, amscd}
\usepackage{mathrsfs}
\usepackage[all]{xy}
\usepackage{amsbsy}
\usepackage{bm}
\usepackage[left=1.2in,top=1.2in,bottom=1.2in,right=1.2in]{geometry}
\usepackage{graphicx}
\usepackage{epstopdf}
\usepackage{verbatim}
\usepackage{multirow}
\usepackage{longtable}

\numberwithin{equation}{section}

\newcommand{\mb}{\mathbf}
\newcommand{\eg}{\emph{e.g.}}
\newcommand{\ie}{\emph{i.e.}}
\newcommand{\cf}{\emph{cf.}}
\newcommand{\ol}{\overline}

\newcommand{\be}{\begin{equation}}
\newcommand{\ee}{\end{equation}}
\newcommand{\ben}{\begin{equation*}}
\newcommand{\een}{\end{equation*}}
\newcommand{\bea}{\begin{eqnarray}}
\newcommand{\eea}{\end{eqnarray}}
\newcommand{\bean}{\begin{eqnarray*}}
\newcommand{\eean}{\end{eqnarray*}}
\newcommand{\nno}{\nonumber}

\newcommand{\bTi}{\begin{itemize} \setlength{\itemsep}{-.1cm}}
\newcommand{\eTi}{\end{itemize}}

\newcommand{\pd}{\partial}
\newcommand{\Tr}{{\rm Tr}}

\renewcommand{\Im}{{\rm Im}}
\renewcommand{\Re}{{\rm Re}}

\newcommand{\Vol}{{\rm Vol}}

\newcommand{\CS}{{\rm CS}}

\newcommand{\Li}{{\rm Li}}

\newcommand{\overcrossing}{{\,\raisebox{-.13cm}{\includegraphics[width=.5cm]{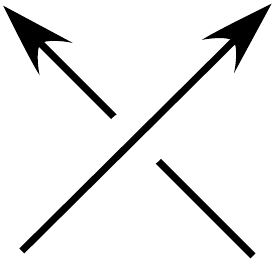}}\,}}
\newcommand{\undercrossing}{{\,\raisebox{-.13cm}{\includegraphics[width=.5cm]{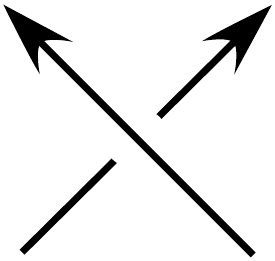}}\,}}
\newcommand{\smoothing}{{\,\raisebox{-.13cm}{\includegraphics[width=.5cm]{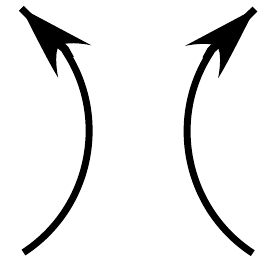}}\,}}
\newcommand{\unknot}{{\,\raisebox{-.08cm}{\includegraphics[width=.4cm]{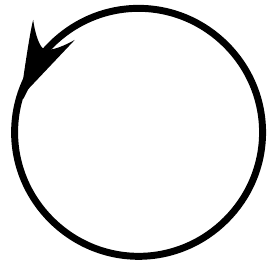}}\,}}

\newcommand{\bs}{\backslash}

\newcommand{\R}{{\mathbb{R}}}
\newcommand{\C}{{\mathbb{C}}}
\newcommand{\Z}{{\mathbb{Z}}}

\newcommand{\Q}{{\mathbb{Q}}}

\renewcommand{\P}{{\mathbb{P}}}

\newcommand{\fg}{{\frak{g}}}

\newcommand{\ft}{{\frak{t}}}

\newcommand{\diag}{{\rm diag}}

\newcommand{\rank}{{\rm rank}}

\newcommand{\Hol}{{\rm Hol}}

\newcommand{\CA}{{\cal A}}
\newcommand{\CB}{{\cal B}}

\renewcommand{\CD}{{\cal D}}

\newcommand{\CH}{{\cal H}}

\newcommand{\CL}{{\cal L}}
\newcommand{\CM}{{\cal M}}

\newcommand{\CO}{{\cal O}}
\newcommand{\CP}{{\cal P}}

\newcommand{\CW}{{\cal W}}

\renewcommand{\a}{\alpha}
\newcommand{\p}{\partial}

\newcommand\bes{\begin{equation*}}
\newcommand\ees{\end{equation*}}
\newcommand{\beas}{\begin{eqnarray*}}
\newcommand{\eeas}{\end{eqnarray*}}

\newcommand{\bse}{\begin{subequations}}
\newcommand{\ese}{\end{subequations}}

\long\def\symbolfootnote[#1]#2{\begingroup%
\def\thefootnote{\fnsymbol{footnote}}\footnote[#1]{#2}\endgroup}

\begin{document}

\pagestyle{empty}

\begin{flushright} 
\end{flushright}
\vspace{.8in}

\begin{center}{\LARGE Quantum Field Theory and the Volume Conjecture} \end{center}

\vspace{.1cm}

\begin{center}{\large Tudor Dimofte and Sergei Gukov} \end{center}
\begin{center} {\it California Institute of Technology 452-48,
Pasadena, CA 91125} \end{center}

\vspace{1in}

The volume conjecture
states that for a hyperbolic knot $K$ in the three-sphere $S^3$
the asymptotic growth of the colored Jones polynomial of $K$
is governed by the hyperbolic volume of the knot complement $S^3\bs K$.
The conjecture relates two topological invariants, one combinatorial
and one geometric, in a very nonobvious, nontrivial manner.
The goal of the present lectures\symbolfootnote[1]{These notes are based on lectures given by the
authors at the workshops \emph{Interactions Between Hyperbolic
Geometry, Quantum Topology, and Knot Theory} (Columbia University, June 2009), \emph{Chern-Simons Gauge Theory: 20 years after} (Hausdorff Center for Mathematics, August 2009), and \emph{Low Dimensional Topology and Number Theory II} (University of Tokyo, March 2010).}
is to review the original statement
of the volume conjecture and its recent extensions and generalizations,
and to show how, in the most general context, the conjecture can be understood
in terms of topological quantum field theory.
In particular, we consider:
$a)$ generalization of the volume conjecture to families of incomplete hyperbolic metrics;
$b)$ generalization that involves not only the leading (volume) term,
but the entire asymptotic expansion in $1/N$;
$c)$ generalization to quantum group invariants for groups of higher rank; and
$d)$ generalization to arbitrary links in arbitrary three-manifolds.

\newpage
\setcounter{page}{1}
\pagestyle{plain}

\noindent\rule{\textwidth}{0.01cm}
\tableofcontents
\noindent\rule{\textwidth}{0.01cm}
\vspace{.5in}

\section{Preliminaries}
\label{sec:basics}

Let $K$ be an oriented knot (or link) in the three-sphere $S^3$.
The original volume conjecture \cite{Kashaev-1997,Mur-Mur}
relates the $N$-colored Jones polynomial of $K$
to the hyperbolic volume of the knot complement $S^3\bs K$:
\be\begin{array}{c@{\qquad}c@{\qquad}c} \text{$N$-colored Jones poly of $K$} & \longleftrightarrow & \text{hyperbolic volume of $S^3\bs K$} \\
\text{(combinatorial, rep. theory)} & & \text{(geometric)\,.}
\end{array} \label{VCtext}
\ee
We begin by reviewing some of the definitions and ingredients
that enter on the two sides here in order to make this statement more precise,
and to serve as a precursor for its subsequent generalization.

\subsubsection*{Jones polynomials}

The (non-colored) Jones polynomial $J(K;q)$ of a knot or link can be defined combinatorially via the skein relation
\be q\,J(\overcrossing) - q^{-1}J(\undercrossing) = (q^\frac12-q^{-\frac12})\,J(\smoothing)\,, \label{skein} \ee
along with the normalization\footnote{The most common normalization for the unknot seen in the mathematics literature is $J(\raisebox{-.05cm}{\includegraphics[width=.3cm]{unknot}})=1$. For the connection with topological quantum field theory, however, \eqref{JU} is much more natural.}
\be J(\unknot) = q^\frac12+q^{-\frac12} \quad\qquad \text{for $\unknot=$ unknot}\,, \label{JU}  \ee
and the rule
\be J(K_1\sqcup K_2) = J(K_1)\,J(K_2)\, \ee
for any disjoint union of knots or links.
Thus, for example, the (right-handed) trefoil and figure-eight knots have Jones polynomials
\begin{align*} J(\mb{3_1}) &= q^{-\frac12}+q^{-\frac32}+q^{-\frac52}-q^{-\frac92}\,, \\
J(\mb{4_1}) &= q^{\frac52}+q^{-\frac52}\,.
\end{align*}
In general, $J(K,q)$ is a Laurent polynomial, $J(K,q)\in \Z[q^\frac12,q^{-\frac12}]$.

The combinatorial construction of the Jones polynomial is intimately related to representation theory of $SU(2)$ --- or the closely related representation theories of the quantum group $U_q(su(2))$ or the affine Lie algebra $\widehat{su(2)}$. In particular, the classical Jones polynomial above is obtained by ``coloring'' the knot (or link) $K$ in $S^3$ with the 2-dimensional representation of $SU(2)$. More generally, such a knot or link can be colored with any finite-dimensional representation $R$ of $SU(2)$, leading to a colored Jones polynomial $J_R(K,q)$. The \emph{N-colored} Jones polynomial $J_N(K,q)$ takes $R$ to be the irreducible $N$-dimensional representation \cite{Wit-Jones, Resh-Tur, Kir-Mel}. The colored Jones polynomial can again be computed in a purely algebraic/combinatorial manner, by using the two rules
\begin{subequations}\label{NJ}
\be
 J_{\oplus_iR_i}(K;q) = \sum_i J_{R_i}(K;q)\, \vspace{-.5cm}\ee
and
\be  J_R(K^n;q) = J_{R^{\otimes n}}(K;q)\,,\,  \ee
\end{subequations}
together with
\be J_R(K_1\sqcup K_2) = J_R(K_1)J_R(K_2)\,, \ee
and the fact that $J_1(K;q)=J_{R=\oslash}(K;q)\equiv 1$.
The first rule says that if $R$ is reducible, then $J_R$ splits as a sum over irreducible components. The second rule says that the $R$-colored Jones polynomial for the $n$-cabling of a knot (formed by taking $n$ copies of the knot or link, slightly displaced away from one another\footnote{This displacement must be done in a way that produces zero linking number between the various copies.}) is equal to the colored Jones polynomial of the original knot but in representation $R^{\otimes n}$.

For example, from (\ref{NJ}a-b) and the fact that $J_{N=2}(K,q) = J(K,q)$, it is easy to see that
\be J_N(\unknot) = \frac{q^{\frac N2}-q^{-\frac N2}}{q^{\frac12}-q^{-\frac12}}\,. \ee
More generally, for any knot $K$, relations \eqref{NJ} can be used to reduce $J_N(K;q)$ to ordinary Jones polynomials of $K$ and its cablings. We have
\begin{align*}
J_1(K;q) &= 1\,, \\
J_2(K;q) &= J(K;q)\,, \\
J_3(K;q) &= J(K^2;q)-1\,, \\
J_4(K;q) &= J(K^3;q)-2J(K;q)\,, \\
 &\;\;\ldots,
\end{align*}
where the expressions for $J_3$, $J_4$, etc. follow from the rules for decomposing representations of $SU(2)$:
$\mb{2}^{\otimes2} = \mb{1}\oplus\mb{3}$\,,\,
$\mb{2}^{\otimes 3}=(\mb{1}\oplus\mb{3})\otimes \mb{2}=\mb{2}\oplus\mb{2}\oplus\mb{4}$\,, etc.
Since $J(K;q) \in \Z[q^\frac12,q^{-\frac12}]$ for any $K$, it is clear that the colored Jones polynomials $J_N(K;q)$ will also be elements of $\Z[q^\frac12,q^{-\frac12}]$\,.

We have explained the left side of \eqref{VCtext}, completely, if somewhat tersely, in terms of algebra and
combinatorics. The right side has a very different interpretation.

\subsubsection*{Hyperbolic volumes}

It was conjectured by Thurston \cite{thurston-1982} (and is now proved \cite{Perelman-geometrization}) that every three-manifold may be decomposed into pieces that admit exactly one of eight different geometric structures. The most common structure by far is hyperbolic. Indeed, in the case of knot complements in $S^3$ this statement can be made exact: a knot complement has a hyperbolic structure if and only if it is not a torus or satellite knot~\cite{thurston-1982}. By definition, a ``hyperbolic structure'' refers to a geodesically-complete metric of constant curvature $-1$. If a hyperbolic structure does exist on a manifold $M$, then it is unique, and the corresponding hyperbolic volume $\Vol(M)$ is a well-defined topological invariant.

In fact, there also exists a natural complexification of the hyperbolic volume of a three-manifold $M$, obtained as
\be \Vol(M) + i \CS(M)\,, \label{VolCS} \ee
where $\CS(M)$ is the so-called Chern-Simons invariant of $M$. To understand this, suppose that $M$ allows a spin structure (as all knot/link complements in $S^3$ do) and consider flat $SL(2,\C)$ connections on $M$ in place of hyperbolic metrics.%
\footnote{Recall that a ``$G$-connection'' on a principal $G$-bundle $E\to M$ can be written locally as a $\fg$-valued one-form $\CA$. The bundle $E$ is typically taken to be trivial in the present context, $E=G\times M$. A gauge transformation (a change of coordinates on $E$) induced by an element $g\in\Gamma(E)$ acts locally on the connection as $\CA\mapsto g^{-1}\CA g+g^{-1}dg$.} %
There exists a flat connection $\CA$ whose real and imaginary parts can be interpreted, respectively, as the vielbein and spin connection of the hyperbolic metric. The real part of the quantity
\be \frac{i}{2}I_{CS}(\CA) = \frac{i}{2}\int_M\Big(\CA\wedge d\CA+\frac{2}{3}\CA\wedge\CA\wedge\CA\Big) \label{ICS} \ee
then reproduces $\Vol(M)$, while the imaginary part defines $\CS(M)$.
The expression $I_{CS}(\CA)$ is the so-called Chern-Simons functional of $\CA$.
Further details can be found \eg\ in \cite{Wit-CSGrav, witten-1991} or \cite{gukov-2003, DGLZ}.
Under gauge transformations acting on $\CA$, the functional $I_{CS}(\CA)$ is only well-defined
up to shifts of $8\pi^2$, leading to an ambiguity of $4\pi^2$ in the definition of $\CS(M)$.
Because of this, it is often convenient to exponentiate the complexified volume \eqref{VolCS},
writing it in the unambiguous form
\be Z(M) =  e^{\frac{i}{4\pi}I_{CS}(\CA)} = e^{\frac{1}{2\pi}\big(\Vol(M)+i\CS(M)\big)}\,. \label{Zcl} \ee

For hyperbolic knot complements, the full complexified volume $Z(M)$ can be efficiently computed in terms of ideal hyperbolic triangulations, \cf\ \cite{snap, neumann-2004, zickert-2008}.

\subsubsection*{The Volume Conjecture}

We have not said much yet about the variable $q$ appearing in the Jones polynomials. Strictly speaking, this variable should be a root of unity%
\footnote{In terms of representation theory, the integer $k$ is identified as the level of the affine Lie algebra $\widehat{su(2)}_k$. The representation theory of the quantum group $U_q(su(2))$ also (crucially) simplifies greatly when $q$ is a root of unity, becoming essentially equivalent to the representation theory of $\widehat{su(2)}_k$. See also Sections \ref{sec:TQFT_CS}-\ref{sec:path}.}
\be q = e^{\frac{2\pi i}{k}}\,, \qquad k \in \Z_+\,. \ee
At the special value $k=N$, all Jones polynomials $J_N(K;q)$ vanish, but the ratio
\be V_N(K,q) =  \frac{J_N(K;q)}{J_N(\unknot;q)} \ee
remains finite. The original volume conjecture \cite{Kashaev-1997, Mur-Mur} then states that
\be \lim_{N\to\infty} \frac{2\pi\log\big|V_N(K;q=e^{\frac{2\pi i}{N}})\big|}{N} = \Vol(M)\,. \ee
It is also possible to remove the absolute value and exponentiate to obtain the complexified generalization (\cf\ \cite{Mur-cx})
\be
\boxed{V_N(K;q=e^{\frac{2\pi i}{N}}) \;\overset{N\to\infty}{\sim}\;  Z(M)^N = e^{\frac{N}{2\pi}\big(\Vol(M)+i\CS(M)\big)}}\,.  \label{VC} \ee

\begin{figure}[ht]
\centering
\includegraphics[width=1.5in]{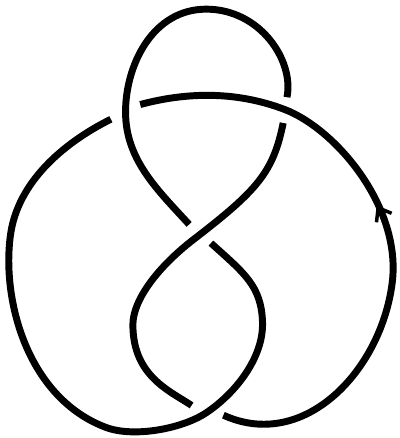}
\caption{The figure-eight knot, $\mb{4_1}$}
\label{fig:41}
\end{figure}

As an example, consider the figure-eight knot (Figure \ref{fig:41}), the simplest hyperbolic knot. The colored Jones polynomial (see \eg\ \cite{Kashaev-1997} or \cite{Habiro-Jones}) is
\be V_N(\mb{4_1};q=e^{\frac{2\pi i}{N}}) = \sum_{m=0}^{N-1}(q)_m(q^{-1})_m\,,\qquad (x)_m := (1-x)(1-x^2)\cdots(1-x^m)\,. \ee
The hyperbolic volume of the figure-eight knot complement is
\be \Vol(S^3\bs\mb{4_1}) = 2\,\Vol(\Delta) = 2.02988...\,, \label{vol41} \ee
where $\Vol(\Delta) = \Im\,\Li_2(e^{i\frac{\pi}{3}})$ denotes the volume of a regular hyperbolic ideal tetrahedron. The Chern-Simons invariant $\CS(S^3\bs\mb{4_1})$ vanishes. It is fairly straightforward (and an informative exercise\footnote{One method involves analytically continuing the summand as a ratio of quantum dilogarithm functions (\cf\ \cite{Fad-Kash, DGLZ}), approximating the sum by an integral, and evaluating it at its saddle point.}) to show that in the limit $N\to\infty$ one has, as expected,
\be \lim_{N\to\infty} \frac{2\pi \log V_N(\mb{4_1};e^{\frac{2\pi i}{N}})}{N} = \Vol(S^3\bs\mb{4_1})\,. \ee

%%%%%%%%%%%%%%%%%%%%%%%%%%%%%%%%%%%%%%%%%%%%%%%%%%%%%%%%%%%%%%%%%%%%%%%%%%%%%%%%%%%%%%%%%%%%%%%%%%%%%

\section{The many dimensions of the volume conjecture}
\label{sec:gen}

There are several natural ways in which one might try to generalize the basic volume conjecture \eqref{VC}.
One possibility is to consider not just $k=N$ (or $q=e^{\frac{2\pi i}{N}}$), but arbitrary values of $k$ (or $q$).
Another option would be to ask what happens to subleading terms in the asymptotic expansion of $V_N(K;q)$ as $N\to\infty$.
It might also be interesting to consider not just hyperbolic knots in $S^3$
but arbitrary links in more complicated three-manifolds.
It turns out that all these generalizations make sense,
and can be nicely combined and interpreted in terms of Chern-Simons theory with complex gauge group \cite{gukov-2003}.
In this section, we detail each of them (and one additional generalization) in turn,
and begin to explain what kind of new objects one should expect on the right-hand-side of \eqref{VC}.
Then, in section \ref{sec:TQFT}, our goal will be to explain where such generalizations come from.

\subsection{Parametrized VC}

The original volume conjecture only held for a special root of unity $q = e^{\frac{2\pi i}{N}}$.
In order to generalize to arbitrary $q=e^{\frac{2\pi i}{k}}$, the appropriate limit to consider is
\be k\to\infty\,,\qquad N\to\infty\,,\qquad u := i\pi\frac { N} k \quad\text{fixed}\, \label{lim} \ee
(or $q\to 1$,\; $q^N = e^{2u}$ fixed). The question, then, is how to understand
\be \lim_{k,N\to\infty} J_N(K;q)^{1/k}\quad ? \ee

\begin{figure}[htb]
\centering
\includegraphics[width=6in]{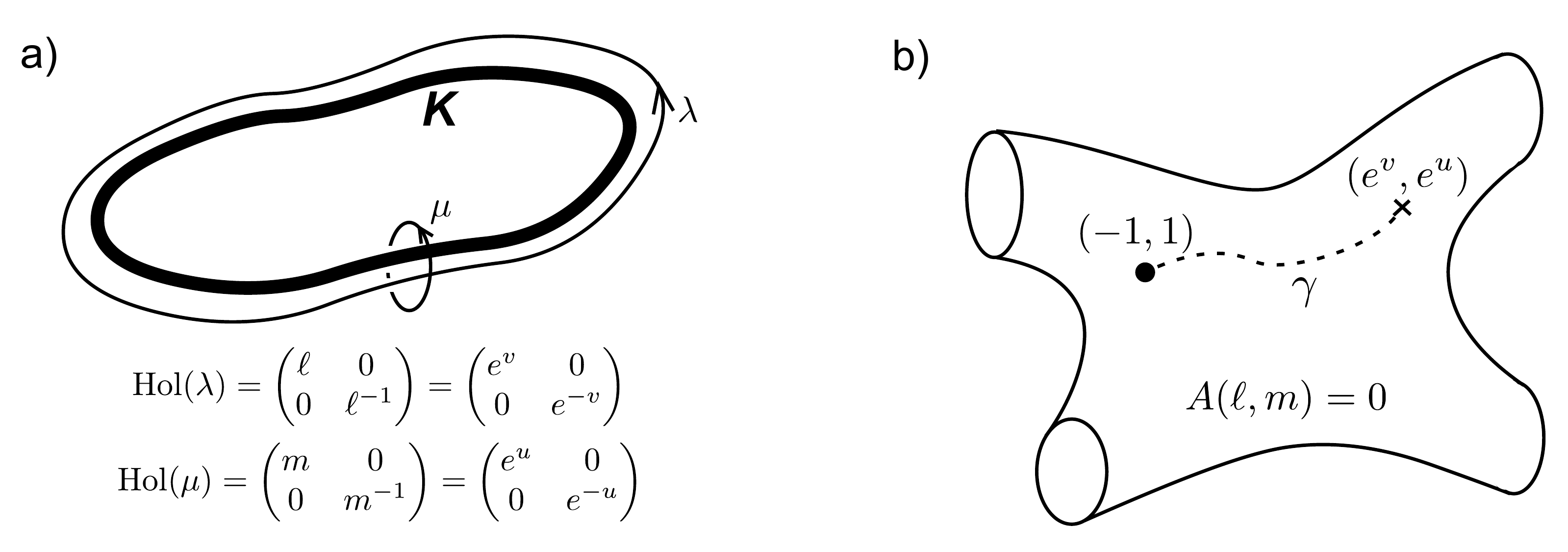}
\caption{ a) The ``longitude'' $\lambda$ and ``meridian'' $\mu$ holonomy paths in the knot complement $S^3\bs K$.  b) Integration on the A-polynomial curve to find the deformed complex volume.}
\label{fig:lm}
\end{figure}

The answer, described in \cite{gukov-2003}, uses the fact that in correspondence
with the ``deformation'' in the colored Jones polynomial, there exists a one-parameter
deformation of the hyperbolic structure on a knot complement $S^3\bs K$.
To understand this, let $\mu$ be a small loop linking the excised knot $K$, as in Figure \ref{fig:lm}a.
In terms of flat $SL(2,\C)$ connections, the geodesically complete hyperbolic metric has a parabolic $SL(2,\C)$ holonomy around $\mu$,
\be \Hol(\mu,{\rm complete}) =  \pm\begin{pmatrix} 1 & 1 \\ 0 & 1\end{pmatrix}\,, \label{hol0} \ee
whereas the incomplete, $u$-deformed hyperbolic metric/$SL(2,\C)$ connection is defined to have a holonomy conjugate to
\be \Hol(\mu,u)  = \begin{pmatrix} e^u & 1 \\ 0 & e^{-u} \end{pmatrix}\,. \label{holu} \ee
(As long as $e^{u}\neq e^{-u}$, this deformed holonomy is also conjugate to the purely diagonal matrix $\diag(e^u,e^{-u})$.)
The resulting metric is not complete.
For example, when $u$ is purely imaginary, the $u$-deformed metric has a conical cusp of angle $2\Im(u)$ at the knot $K$.

The complexified hyperbolic volume for this one-parameter family of metrics can again be defined in terms of the Chern-Simons functional $I_{CS}(\CA)$ appearing in \eqref{ICS}.
Now, however, $\CA=\CA(u)$ should be a flat $SL(2,\C)$ connection with prescribed holonomy \eqref{holu}.
The "parametrized" volume conjecture then takes the form \cite{gukov-2003}
\be \boxed{ J_N(K;q) \;\overset{k,N\to\infty}{\sim}\; e^{-\frac{k}{4\pi i}I_{CS}(\CA(u))} } \,. \label{VCu} \ee

The quantity $I_{CS}(\CA(u))$ can be described very explicitly. Indeed, suppose that we require a hyperbolic metric (expressed in terms of a flat $SL(2,\C)$ connection) to have holonomies conjugate to $\diag(e^u,e^{-u})$ and $\diag(e^v,e^{-v})$, respectively, along the meridian and longitude loops depicted in Figure \ref{fig:lm}a. Such a metric exists if and only if the so-called A-polynomial of $K$ vanishes \cite{cooper-1994},
\be A(\ell,m) = 0 \qquad\text{for}\qquad \ell=e^v\,,\quad m=e^u\,. \label{A} \ee
Given a fixed $e^{u}\in\C^*$, exactly one of the solutions $v=v^{\rm hyp}(u)$ of this equation corresponds to the $u$-deformed hyperbolic metric.
The Chern-Simons functional evaluated at the flat connection $\CA(u)$ can then be written as
\be I_{CS}(\CA(u)) = I_{CS}(\CA(i\pi)) +4 \int_\gamma \theta\,, \label{CStheta} \ee
where $\CA(i\pi)$ is the non-deformed hyperbolic flat connection,
\be \theta = -(v'+i\pi)\,du' \ee
is a one-form on the curve $A(e^{v'},e^{u'})=0$, and $\gamma$ is a path on this curve that connects the complete hyperbolic structure at $(e^{v'},e^{u'})=(-1,\pm1)$ to the $u$-deformed metric at $(e^{v'},e^{u'})=(e^{v^{\rm hyp}(u)},e^u)$, as in Figure \ref{fig:lm}b.%
\footnote{The actual complexified volume that appears in the literature on hyperbolic geometry (\cf\ \cite{NZ, yoshida-1985, Hilden-1996}) is related to $I_{CS}(\CA(u))$ as
\be \Vol(S^3\bs K;u)+i\CS(S^3\bs K;u) = \frac{i}{2}I_{CS}(\CA(u)) + 2iv(u)\Re(u)-2\pi u+2\pi^2 i\,.
\ee
Note that $I_{CS}(\CA(u))$ is analytic in $u$, whereas $\Vol(u)+i\CS(u)$ is not.}

As our recurrent example, consider again the figure-eight knot. The complete colored Jones polynomial, \cf\ \cite{Habiro-Jones}, is
\be J_N(\mb{4_1};q) = \frac{q^{\frac N 2}-q^{-\frac N 2}}{q^\frac12-q^{-\frac12}}\sum_{j=0}^{N-1}q^{Nj}\prod_{k=1}^j(1-q^{k-N})(1-q^{-k-N})\,. \ee
The A-polynomial of the figure-eight knot is
\be A(\ell,m) = (\ell-1)(m^4\ell^2-(1-m^2-2m^4-m^6+m^8)\ell+m^4\ell^2)\,,\ee
and from \eqref{CStheta} and \eqref{vol41}, it results (after some algebra) that the Chern-Simons functional can be written as
\be I_{CS}(\CA(u)) = 2\Li_2(e^{-p-u})-2\Li_2(e^{p-u})+8(p-i\pi)(u-i\pi)\,, \ee
where $x=e^p$ is the solution to $m^3x^2+(1-m^2-m^4)x+m^3=0$ with smallest negative imaginary part.
For irrational $u/ i\pi$ in a neighborhood of $u=i\pi$ it can then be shown (\cf\ \cite{Mur-parametrized, Mur-Yok}) that the proposed asymptotics \eqref{VCu} indeed hold. \\

The necessity for taking $u/i\pi$ irrational here may appear a little strange at first glance.
It stems fundamentally from the fact that the Jones polynomials $J_N(K;q=e^{\frac{2\pi i}{k}})$
are really only defined for $N,k\in \Z$.
A subtle analytic continuation in either $N$ or $k$ is necessary to achieve $u/i\pi = N/k\notin \Q$.
As anticipated in \cite{gukov-2003} and explained recently in \cite{Wit-anal},
it is this continuation that causes the growth of the colored Jones polynomial to be exponential.
We will remark on this further in Section \ref{sec:path}.

In light of this argument, one might ask now why the original volume conjecture
at the rational value $k=N$ or $u=i\pi$ held in the first place.
Recall that $J_N(K;q)$ actually vanished at $k=N$, so it was necessary
to divide by $J_N(\unknot;q)$ to obtain the non-vanishing ratio $V_N(K;q)$.
Examining $V_N(K;q)$ at $u\to 0$ is equivalent to considering
the \emph{derivative} of $J_N(K;q)$ at $u=i\pi$, which of course knows about
analytic continuation.\footnote{We thank E. Witten for useful observations on this subject.}

\subsection{Quantum VC}

The second option for generalizing the volume conjecture \eqref{VC} is to ask for higher-order terms in the asymptotic expansion of the colored Jones polynomial. Let us define a new ``quantum'' parameter $\hbar$ as
\be  \hbar = \frac{i\pi}{k}\,, \ee
so that
\be q = e^{2\hbar}\,. \ee
The two parameters $N$ and $k$ of the colored Jones polynomial can be traded for $\hbar$ and $u$, and the limit \eqref{lim} is simply $\hbar\to 0$. At $u=i\pi$, higher-order asymptotics are then predicted \cite{gukov-2003, DGLZ} to have the form
\be V_N(K;q=e^{\frac{2\pi i}{N}}) \;\overset{N\to\infty}{\sim}\;
 \exp\left(\frac{1}{2\hbar}(\Vol+i\CS)-\frac{3}{2}\log\hbar+\frac{1}{2}\log\frac{-i\pi T_K}{4} + \sum_{n=2}^\infty \tilde{S}_n\hbar^{n-1}\right)\,. \label{VCh0}
\ee
Here, for example, $T_K$ is the Ray-Singer torsion of the knot complement $S^3\bs K$. It can be defined after putting any background metric on $S^3\bs K$ \cite{Porti} as
\be T(M) = \exp\left(-\frac12\sum_{n=0}^3n(-1)^n\log{\rm det}'\Delta_n\right) = \frac{(\det'\Delta_0)^\frac32}{(\det'\Delta_1)^\frac12}\,,\ee
where $\Delta_n$ is the Laplacian acting on $n$-forms.

It is fairly straightforward to combine the present quantum deformation with the parametrization of the volume conjecture in $u$. The expectation is that
\be \boxed{J_N(K;q) \;\overset{N,k\to\infty}{\sim}\;
 \exp\left( -\frac{1}{4\hbar}I_{CS}(\CA(u)) -\frac{3}{2}\log\hbar + \frac{1}{2}\log \frac{iT_K(u)}{4\pi}+\sum_{n=2}^\infty S_n(u)\hbar^{n-1}\right)}\,. \label{VCh}
\ee
Here, $T_K(u)$ is a $u$-deformed torsion, and is related to the Alexander polynomial of $K$ \cite{Mur-Alex}.
The higher-order coefficients in \eqref{VCh0} are related to those in \eqref{VCh} as
\be \sum_{n\geq2}\tilde{S}_n\hbar^{n-1}=\sum_{n\geq2}S_n(i\pi)\hbar^{n-1}-\log\frac{\sinh\hbar}{\hbar}\,. \ee

For the figure-eight knot, the quantum volume conjecture \eqref{VCh} was tested to first subleading order in \cite{gukov-2006}, using the Ray-Singer torsion
\be T_{\mb{4_1}}(u) = \frac{4\pi^2}{\sqrt{-m^{-4}+2m^{-2}+1+2m^{2}-m^{4}}}\,. \ee
Higher-order coefficients $S_n(u)$ can also be computed \cite{DGLZ}. For example,
\begin{align} S_2(u) &= \frac{-i(T_{\mb{4_1}})^3}{12(4\pi^2)^3m^6}\big(1-m^2-2m^4+15m^6-2m^8-m^{10}+m^{12}\big)\,, \\
  S_3(u) &= \frac{-2(T_{\mb{4_1}})^6}{(4\pi^2)^6m^6}\big(1-m^2-2m^4+5m^6-2m^8-m^{10}+m^{12}\big)-\frac16\,.
\end{align}
These expressions appear to be new, unexplored knot invariants with distinctive number-theoretic properties.
Needless to say, it would be interesting to test the quantum volume conjecture \eqref{VCh}
for other hyperbolic knots and/or to higher order in the $\hbar$-expansion. \\

Just as the generalization of the volume conjecture to $u\neq 0$
was interpreted in terms of the $SL(2,\C)$ Chern-Simons functional,
there is also a Chern-Simons interpretation of the quantum volume conjecture.
One must consider how the functional $I_{CS}(\CA)$ behaves when
the connection $\CA$ undergoes ``quantum fluctuations'' away from the flat connection $\CA(u)$.
This is accomplished in physics via perturbative quantum field theory.
Symbolically, we can write $\CA=\CA(u)+\CA'$, where $\CA'$ contains
the fluctuations away from flatness, and define a perturbative ``partition function'' via the path integral
\be Z(S^3\bs K;u;\hbar)_{\rm pert} = \int \CD\CA'\,e^{-\frac{1}{4\hbar}I_{CS}(\CA(u)+\CA')}\,. \label{path2} \ee
The exponent in the integrand has a critical point at $\CA'=0$, and a saddle point expansion around this point yields the right-hand-side of \eqref{VCh}.
(To be very precise, $J_N(K;q)\sim Z(S^3\bs K;\hbar;u)/Z(S^3;\hbar)$, where $Z(S^3;\hbar) = \sqrt{2/k}\sin(\pi/k)$ is the partition function of the three-sphere $S^3$.)

%%%%%%%%%%%%%%%%%%%%%%%%%%%%%%%%%%%%%%%%%%%%%%%%%%%%%%%%%%%%%%%%%%%%%%%%%%%%%%%%%%%%%%%%%%%%%%%%%%%%

\subsection{Groups and representations}

\begin{figure}[ht]
\centering
\includegraphics[width=3.5in]{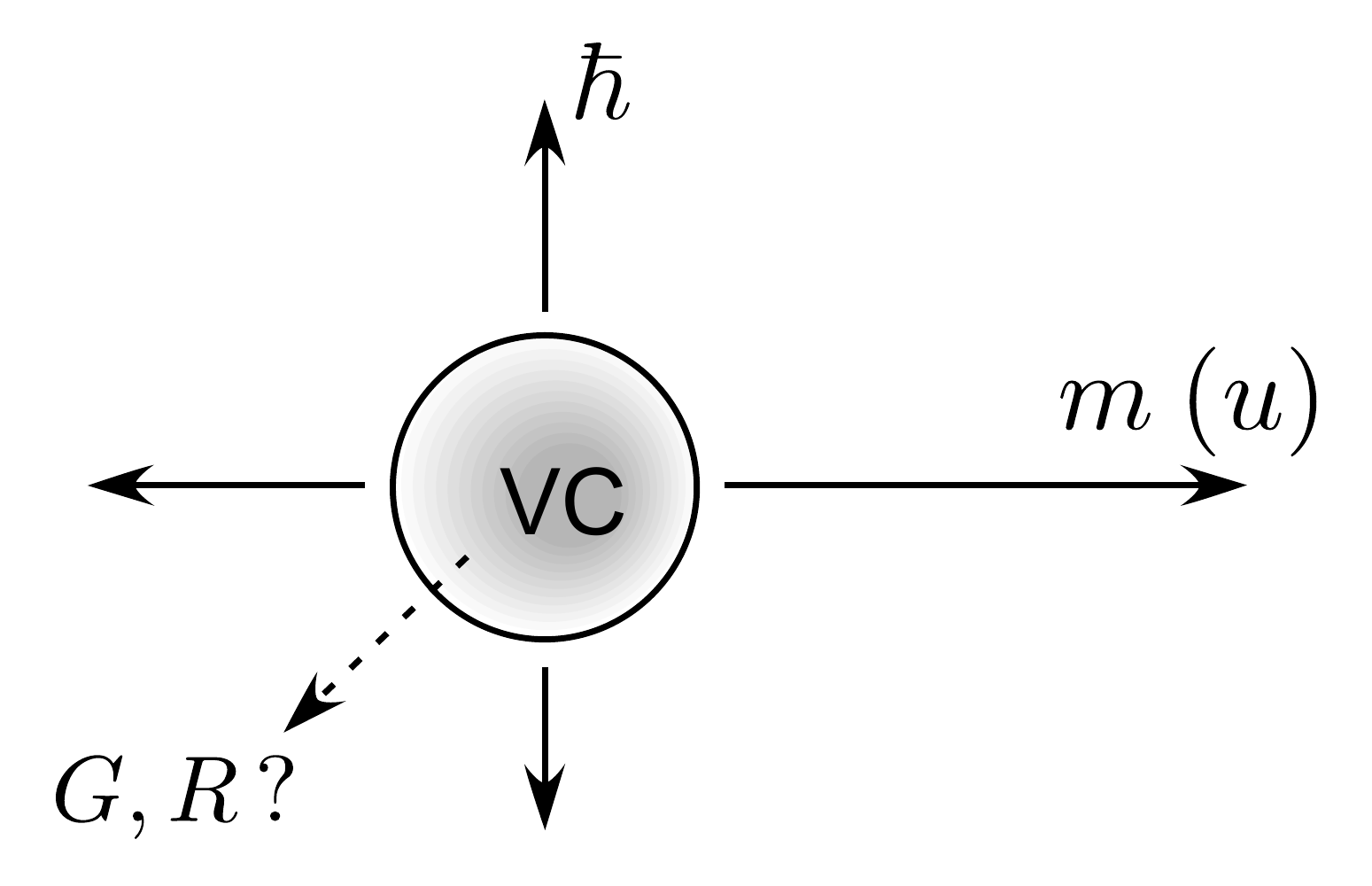}
\caption{Continuous and discrete generalizations of the volume conjecture.}
\label{fig:uh}
\end{figure}

So far, we have considered two continuous deformations of the volume conjecture, in $u$ and $\hbar$, as drawn schematically in Figure \ref{fig:uh}. In addition, there are two discrete generalizations that we can make.

The first such generalization involves the ``gauge groups'' and representations that define colored Jones polynomials. Recall from Section \ref{sec:basics} that the $N$-colored Jones polynomial is a quantum $SU(2)$ invariant that corresponds to coloring a knot with the $N$-dimensional representation of $SU(2)$.
More generally, one can consider ``quantum $SU(n)$ invariants,'' or in fact invariants for any compact Lie group $G$. Knots or links should then be colored by finite-dimensional representations $R$ of $G$. For semisimple $G$ and irreducible $R$, the representation can be labelled by a highest weight $\lambda$ in the weight lattice $ \Lambda_{wt}\subset \fg^*$, where $\fg = {\rm Lie}(G)$. The resulting quantum polynomial invariant of a knot in $S^3$ may be denoted
\be P^G_{R_\lambda}(K;q)\,. \ee

Just like the colored Jones polynomial, $P^G_{R}(K;q)$ depends on a root of unity $q=e^{\frac{2\pi i}{k}}$. Also like the colored Jones, these invariants satisfy
\be  P^G_{\oplus_i R_i}(K;q) = \sum_i P^G_{R_i}(K;q)\,,\qquad\text{and}\qquad P^G_{R^{\otimes n}}(K;q) = P^G_R(K^n;q)\,.\qquad \ee
More general tensor products can also be produced by cabling a knot or link and coloring each component of the cable with a different representation. When $G=SU(n)$ and $R$ is the fundamental representation (or any of its conjugates), the polynomial $P^G_{R}(K;q)$ satisfies a skein relation similar to \eqref{skein}.

Using the positive nondegenerate trace form $-\Tr\,:\,\fg\times\fg\to\R$, the weight $\lambda$ can be identified with its dual element $\lambda^*$ in $\ft$, the Cartan subalgebra of $\fg$. Let us also define $\rho$ to be half the sum of positive roots, and $\rho^*\in\ft\subset\fg$ its dual. Then the interesting limit to consider for $P^G_{R_\lambda}(K,q)$ is
\be k\to \infty\,,\qquad \lambda^*\to\infty\,,\qquad u:=\frac{i\pi}{k}(\lambda^*+\rho^*)\;\;\text{fixed}\,, \ee
or
\be q = e^{2\hbar}=e^{\frac{2\pi i}{k}}\to 1\,\quad (\hbar\to 0)\,,\qquad q^{\lambda^*+\rho^*}=e^{2u}\;\;\text{fixed}\,.
\ee
The parameter $u$ has now become a diagonal matrix, an element of $\ft_\C$. Coming back to the case of $SU(2)$ and an $N$-dimensional representation, in this notation we have
\be \lambda^*=\begin{pmatrix} N-1 & 0 \\ 0 & -(N-1) \end{pmatrix}\,,
 \qquad \rho^* = \begin{pmatrix} 1 & 0 \\ 0 & -1 \end{pmatrix}\,,
 \qquad u 
  = i\pi\begin{pmatrix} \frac{N}{k} & 0 \\ 0 & -\frac{N}{k}\end{pmatrix}\,.
\ee

The asymptotics of the invariant $P^G_{R}(K;q)$ should look very similar to those of the colored Jones polynomial, namely
\be \boxed{ P^G_{R_\lambda}(K;q) \;\overset{\hbar\to 0}{\sim}\; \exp\left( -\frac{1}{4\hbar}I_{CS}(u)-\frac{\delta}{2}\log\hbar+\frac12\log \frac{iT(u)}{4\pi}+\sum_{n=2}^\infty S_n(u)\hbar^{n-1} \right) } \,. \label{VCRG}
\ee
The leading term $I_{CS}(u)$ is now the Chern-Simons functional \eqref{ICS} evaluated at a flat $G_\C$ connection $\CA(u)$ --- in other words, a connection taking values in the complexified Lie algebra $\fg_\C$ --- whose holonomy around the meridian of the knot as in Figure \ref{fig:lm}a is
\be \Hol(\mu) = m=e^u\,. \ee
For generic $u$, this holonomy is an element of the complexified maximal torus $T_\C \subset G_\C$. Again, $I_{CS}(u)$ can be expressed as
\be I_{CS}(u) = \text{const.} + 4\int_{\gamma(u)} \theta\,, \ee
where\, $\theta \sim -\sum_{i=1}^r v_i\,du_i$ + exact \,\,is a differential on an $r$-dimensional complex variety cut out by $r$ equations $A_j(e^{v},e^{u})=0$, with $r=\rank(G)$. The equations $A_j(e^{v},e^{u})=0$ describe the moduli space of flat $G_\C$ connections on $S^3\bs K$.

Subleading terms on the right side of \eqref{VCRG} also have a geometric interpretation.
The function $T(u)$ is the Ray--Singer torsion of the knot complement twisted
by the flat connection $\CA(u)$,
and the number $\delta$ is a fixed integer which can be computed in terms of
cohomology of $S^3\bs K$ with coefficients in the associated flat bundle,
with structure group $G_\C$ and connection $\CA(u)$ (\cf\ \cite{barnatan-1991w, DGLZ}).
More generally, the full asymptotic expansion can be written as a perturbative path integral
just like \eqref{path2}, which takes into account the quantum fluctuations of a flat $G_\C$ connection.

%%%%%%%%%%%%%%%%%%%%%%%%%%%%%%%%%%%%%%%%%%%%%%%%%%%%%%%%%%%%%%%%%%%%%%%%%%%%%%%%%%%%%%%%%%%%%%%%%%

\subsection{Links and 3-manifolds}
\label{sec:links}

The final generalization of the volume conjecture that we consider is to arbitrary links in arbitrary three-manifolds. Here we really begin to require a true TQFT description of the ``quantum $G$-invariants'' of knots and links. This was supplied by quantum Chern-Simons theory with compact gauge group $G$ in \cite{Wit-Jones}, and reinterpreted via quantum groups and R-matrices in \cite{Resh-Tur}. Using either of these approaches, one may define a quantum partition function
\be  Z^G(M,L;\{R_a\};\hbar) \label{ZM} \ee
for a link $L$ in any three-manifold $M$, where each component of the link is colored with a different representation $R_a$. The ``polynomial'' $P^G_R$ is obtained from this after normalizing by the partition function of an empty manifold,
\be P^G_{\{R_a\}}(M,L;q) = \frac{Z^G(M,L;\{R_a\};\hbar)}{Z^G(M;\hbar)}\,,\qquad (q=e^{2\hbar})\,. \ee
Thus, in the case of the colored Jones polynomial,
\be J_N(K;q) = \frac{Z^{SU(2)}(S^3,K;R_N;\hbar)}{Z^{SU(2)}(S^3;\hbar)}\,. \ee
The integer $k$ (appearing in $q=e^{2\hbar}=e^{\frac{2\pi i}{k}}$) is identified with the ``level'' or coupling constant of the compact Chern-Simons theory.\footnote{To be completely precise, the integer $k$ used throughout these lectures is the sum of the Chern-Simons level and the dual Coxeter number of $G$. \label{foot:k}}

The partition function \eqref{ZM} supplies the left-hand-side of the volume conjecture.
We then want to understand the asymptotics of $Z^G(M,L;\{R_a\};\hbar)$ in the limit $\hbar\to 0$,
with a parameter $u_a = \hbar(\lambda_a^* + \rho^*)$ held fixed for each separate link component.
The answer should be given by perturbative, quantum Chern-Simons theory with \emph{complex} gauge group $G_\C$,
evaluated on the link complement $M\bs L$, in the background of a flat connection with fixed holonomy
\be m_a = \exp(u_a) \ee
around the meridian of each excised link component. Denoting this perturbative Chern-Simons partition function by
\be Z^{G_\C}_{\rm pert}(M\bs L;\{u_a\};\hbar) = \exp\left(-\frac1{4\hbar}I_{CS}(\{u_a\})-\frac\delta2\log\hbar+\ldots\right)\,, \ee
we expect that
\be \boxed{Z^G(M,L;\{R_a\};\hbar) \;\overset{\hbar\to 0}{\sim}\; Z^{G_\C}_{\rm pert}(M\bs L;\{u_a\};\hbar)}\,. \label{fullVC} \ee

This discussion can also be rephrased in a somewhat more symmetric manner, using link complements on both sides of the volume conjecture. It turns out that in compact Chern-Simons theory the partition function of a knot (or link) $K\subset M$ colored by representation $R_\lambda$ is equivalent to the partition function of the knot complement $M\bs K$ with fixed meridian holonomy
\be m = \exp\left(i\pi\frac{\lambda^*+\rho^*}{k}\right) = \exp\left(\hbar(\lambda^*+\rho^*)\right) = \exp(u)\,. \ee
For the compact $G$ theory to make sense, the eigenvalues of the matrix $u/i\pi$ must clearly be rational. However, interesting asymptotics --- potentially with exponential growth as in \eqref{fullVC} --- occur when $u$ is \emph{analytically continued} away from such rational values. This process of analytic continuation naturally lands one in the regime of Chern-Simons theory with complex gauge group $G_\C$~\cite{DGLZ}. \\

After so many generalizations,
it may be unclear that the volume conjecture has anything to do with volumes anymore.
Indeed, for higher-rank gauge groups $G$,
``volume'' should not be a hyperbolic volume but rather
the ``volume'' of a holonomy representation
\be \varrho : \pi_1 (M\bs K) \to G_{\C}\,. \label{rhofirst} \ee
Even in the case of $G=SU(2)$ and knots in the three-sphere,
one may run across cases of non-hyperbolic knot complements.
It was clear from the initial days of the volume conjecture~\cite{Mur-Mur}
that even in these cases the asymptotics of $J_N(K;q)$
could still be given by an appropriate flat (but non-hyperbolic/non-metric) $SL(2,\C)$ structure.

%%%%%%%%%%%%%%%%%%%%%%%%%%%%%%%%%%%%%%%%%%%%%%%%%%%%%%%%%%%%%%%%%%%%%%%%%%%%%%%%%%%%%%%%%%%%%%%%%%%%%%%%

\section{TQFT}
\label{sec:TQFT}

We have just seen that the volume conjecture admits a multitude of generalizations, all of which seem to be related to Chern-Simons quantum field theory. The most complete statement of the volume conjecture \eqref{fullVC} involves Chern-Simons theory with compact gauge group $G$ on the left-hand side and Chern-Simons theory with complex gauge group $G_\C$ on the right:
\be \begin{array}{ccc}
\underline{\text{\emph{combinatorics/rep. theory}}} && \underline{\text{\emph{geometry}}} \vspace{.2cm} \\
\text{quantum $G$-invariants} & & \text{volumes of representations} \\
  J_N(K;q)\,,\;P^G_{R_\lambda}(K;q)\,, &\quad\longleftrightarrow\quad & \varrho : \pi_1 (M\bs K) \to G_{\C} \,, \\
  Z^{G}(M,K;u;\hbar)\,,\;\text{etc.} &&  Z^{G_\C}_{\rm pert}(M\bs K;u;\hbar)\,,\;\text{etc.} \vspace{.4cm}  \\
\multicolumn{3}{c}{q=e^{\frac{2\pi i}{k}}=e^{2\hbar}\,,\quad u = i\pi\frac{\lambda^*+\rho^*}{k}\,. \vspace{.1cm}}
\end{array} \label{VCsum} \ee
Chern-Simons theory is a topological quantum field theory (TQFT).
In addition to the basic implication that partition functions such as $Z^{G}(M,K;u;\hbar)$ or $Z^{G_\C}_{\rm pert}(M\bs K;u;\hbar)$ are topological invariants of colored knots and links in three-manifolds, the structure of TQFT provides powerful methods for actually computing them in multiple ways.
It also shows why a general correspondence like \eqref{VCsum} might be expected to hold.

%%%%%%%%%%%%%%%%%%%%%%%%%%%%%%%%%%%%%%%%%%%%%%%%%%%%%%%%%%%%%%%%%%%%%%%%%%%%%%%%%%%%%%%%%%%%%%%%%%%%%%%%%%%%

\subsection{Cutting and gluing}
\label{sec:cut}

In its more mathematical incarnation, a 3-dimensional TQFT can be thought of as a functor $Z$ that assigns
\be \begin{array}{c@{\qquad}c@{\qquad}c}
 \text{closed 3-manifold $M$} & \leadsto & \text{number $Z(M)$} \\
 \text{closed 2-manifold $\Sigma$} &\leadsto & \text{vector space $Z(\Sigma)$} \\
 \text{closed 1-manifold $S^1$} &\leadsto & \text{category $Z(S^1)$} \\
 \text{point $p$} &\leadsto & \text{2-category $Z(p)$}\,.
\end{array}
\ee
For our applications to Chern-Simons theory, we will really only need the top two levels $Z(M)$ and $Z(\Sigma)$. The finer structure of categories and 2-categories has recently been explored in \eg\ \cite{Freed-TQFT}.

\begin{figure}[htb]
\centering
\includegraphics[width=5.5in]{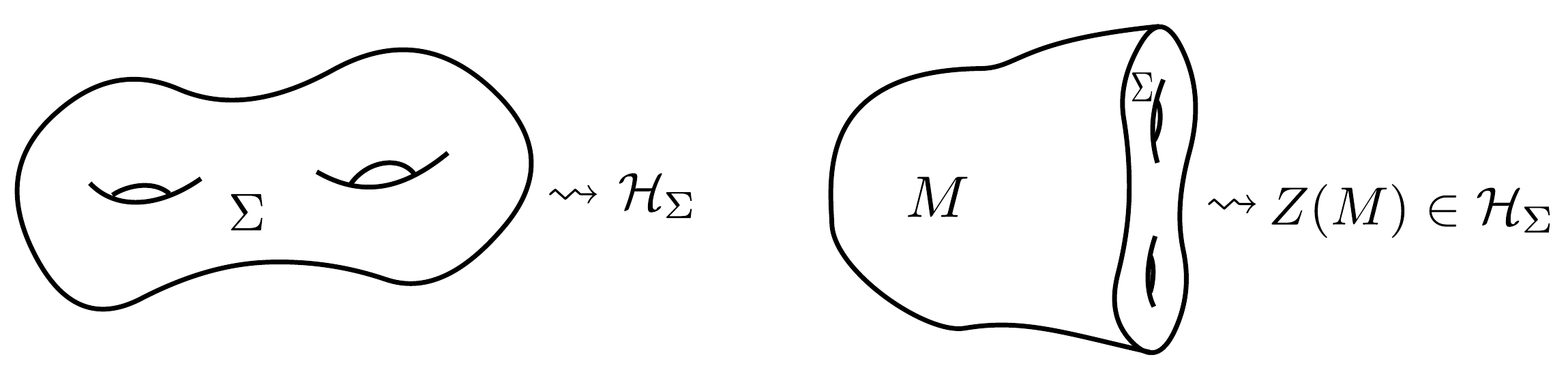}
\caption{Hilbert space assigned to a surface $\Sigma$ and partition function assigned to a three-manifold $M$ in TQFT.}
\label{fig:TQFT}
\end{figure}

If a 3-manifold $M$ has a boundary $\Sigma = \pd M$, the object $Z(M)$ is no longer a number, but an element of the vector space $Z(\Sigma)$ assigned to the boundary, as shown in Figure \ref{fig:TQFT}. This vector space is in fact a Hilbert space, so let us denote it as $\CH_\Sigma = Z(\Sigma)$. At the top two levels, the TQFT must then satisfy the following axioms of Atiyah and Segal (\cf\ \cite{atiyah-1990}).
\begin{enumerate}

\item A change of orientation $\Sigma\to - \Sigma$ dualizes the Hilbert space, $\CH_{-\Sigma}=\CH_\Sigma^*$\,.

\item For a boundary consisting of multiplet disjoint components, $\CH_{\Sigma_1 \sqcup \Sigma_2} = \CH_{\Sigma_1}\otimes\CH_{\Sigma_2}\,.$

\item Using the first two axioms, we see that for a manifold $M$ with $\pd M = (-\Sigma_1)\sqcup\Sigma_2$ one obtains a map $Z(M)\,:\;\CH_{\Sigma_1}\to\CH_{\Sigma_2}$. Then, given a 3-manifold $N$ that can be written as $N = M_1\cup_{\Sigma_2} M_2$, with $\pd M_1 = (-\Sigma_1)\sqcup \Sigma_2$ and $\pd M_2 = (-\Sigma_2)\sqcup \Sigma_3$ as illustrated below, the functoriality property $Z(N) = Z(M_2)\circ Z(M_1)$ must hold. \\
\be \includegraphics[width=3.5in]{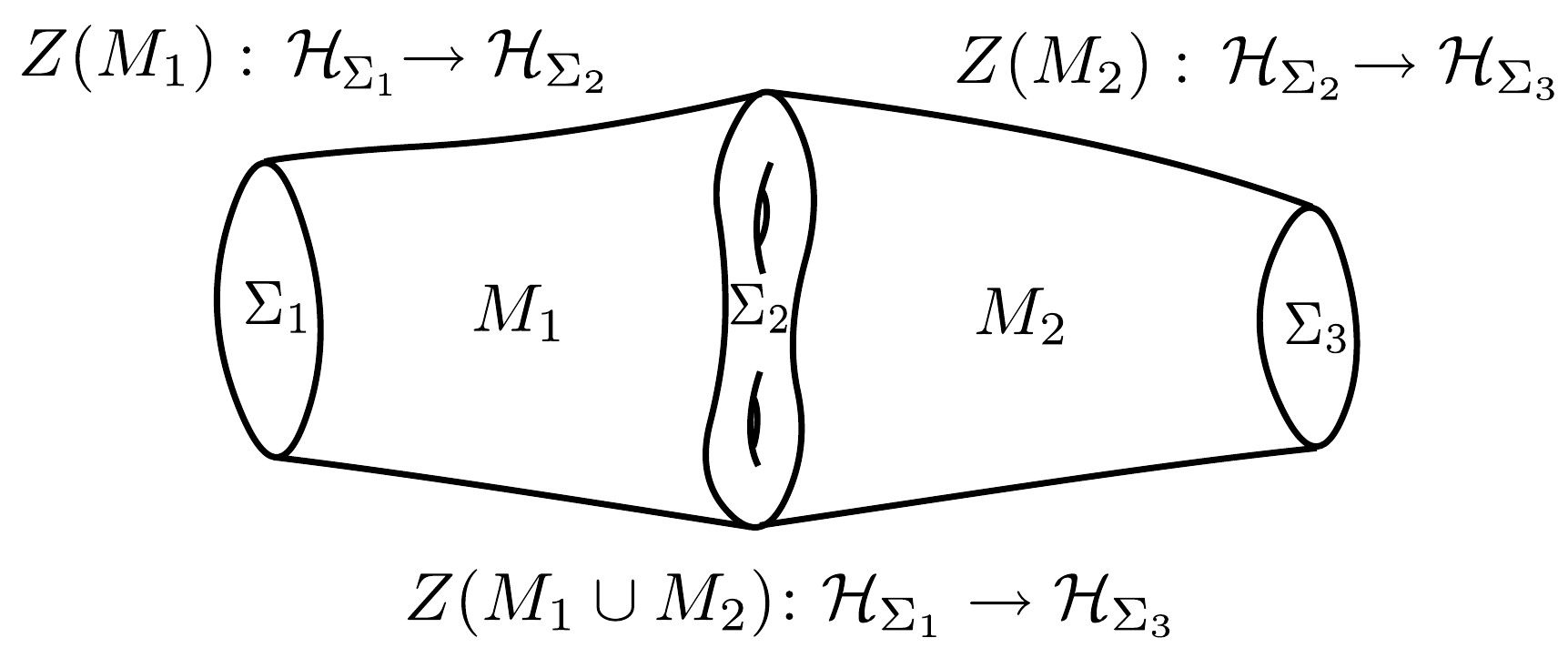} \nno\ee

\item For the empty boundary, $\CH_{\Sigma=\oslash}=\C$\,.

\item For $M=\Sigma\times I$, the map $Z(M)\,:\CH_\Sigma\overset{\rm id}{\to}\CH_\Sigma$ is just the identity.

\end{enumerate}

Using these axioms, the partition function $Z(M)$ of any three-manifold, with or without boundary, may be constructed by cutting the manifold into pieces and taking inner products in boundary Hilbert spaces to glue the pieces back together. For this purpose, it is often convenient to know how the mapping class group of a surface $\Sigma$ acts on $\CH_\Sigma$, in order to properly identify the Hilbert spaces on two sides of a gluing.

There are many examples of three-dimensional TQFT, differing essentially in the definitions of the boundary Hilbert spaces $\CH(\Sigma)$, as well as the action of the mapping class groups on these spaces. In the case of Chern-Simons theory with gauge group $G$ (whether compact or complex), $\CH_\Sigma$ is a quantization of the space $\CM_{{\rm flat}} (G;\Sigma)$ of flat $G$-connections on $\Sigma$:
\be  \CM_{{\rm flat}} (G;\Sigma) = \left\{\begin{array}{r|l}
  \text{connections $\CA$ on principal} & \multirow{2}{*}{$F_\CA=0$} \\
  \text{$G$-bundle over $\Sigma$} &
  \end{array}\right\}\Big/ \text{gauge equivalence}\,.
\ee
(Recall that a connection is flat if the curvature $F_\CA=d\CA+\CA\wedge \CA$ vanishes.)
The precise meaning of the quantization used to obtain $\CH_\Sigma$ from $\CM_{{\rm flat}} (G;\Sigma)$ will be the subject of Section \ref{sec:quant}. It depends on the level $k=i\pi\hbar^{-1}$ (or coupling constant) of Chern-Simons theory, the only adjustable parameter in the TQFT.

In Chern-Simons theory, one is also interested in colored knots or links embedded in 3-manifolds.
Suppose for the moment that we have compact Chern-Simons theory with gauge group $G$ and level $k\in \Z$.
The intersection of a knot and a boundary surface $\Sigma$ shows up as a puncture on $\Sigma$ and TQFT would assign the boundary $S^1$ surrounding this puncture in $\Sigma$ the category of representations of the affine Lie algebra $\widehat{\frak g}_{k}$,
\be Z(S^1) \sim \text{reps of $\widehat{\frak g}_k$}\,. \ee
The definition of the Hilbert space $\CH(\Sigma)$ of a multiple-punctured $\Sigma$ would then have to be altered to include the space of homomorphisms between such representations. For our purposes, however, the complication of knots can be conveniently avoided by \emph{excising the knots and trading representations that color the knots for boundary conditions on knot complements}.

This trick was already mentioned in Section \ref{sec:links}. In the language of TQFT, it can be described the following way. Suppose that we have a knot $K$ colored by representation $R_\lambda$ inside the closed manifold $M$.
We cut out a tubular neighborhood $N_K$ of the knot, so that
\be M = (M\bs N_K)\cup_{T^2}N_K\,,\qquad N_K\simeq D^2\times S^1\,.\ee
Of course, $M\bs N_K\simeq M\bs K$ is just the knot complement, and $N$ is topologically a 2-disk times $S^1$ that contains the knot running through its center. The partition functions $Z(M\bs K;u;\hbar)$ and $Z(N_K;R_\lambda;\hbar)$ are both vectors in the boundary Hilbert space $\CH_{T^2}$; therefore, by TQFT,
\be Z(K\subset M;R_\lambda;\hbar) = \langle\, Z(M\bs K;u;\hbar) \,,\, Z(N_K;R_\lambda;\hbar)\,\rangle_{\CH_{T^2}}\in\C\,.\ee
As we will see in the next section,
the Hilbert space $\CH_{T^2}$ can be understood as a space of functions
of the variable $u$ that describes the holonomy of flat connections around
the meridian of $T^2$ (as in Figure \ref{fig:lm}a).
The crucial fact, then, is that the vector $Z(N;R_\lambda;\hbar)\in \CH_{T^2}$
is only supported on the part of this space with
\be e^u = \Hol(\mu) = \exp\left(i\pi \frac{\lambda^*+\rho^*}{k} \right)\,. \label{holagain} \ee
In other words, $Z(N;R_\lambda;\hbar)$ acts like a delta-function $\delta(u-i\pi \frac{\lambda^*+\rho^*}{k})$. Therefore, coloring by $R_\lambda$ is equivalent to restricting $Z(M\bs K)$ to an appropriate one-dimensional subspace of $\CH_{T^2}$:
\be Z(K\subset M;R_\lambda;\hbar)
= Z(M\bs K;u;\hbar)\big|_{u=i\pi \frac{\lambda^*+\rho^*}{k}} \, \in\, \C \,. \label{uR} \ee \\
Our plan now is to give a complete description of $\CH_{T^2}$
and to explain how the elements $Z(M\bs K) \in \CH_{T^2}$ may be calculated for knot complements,
in the case of Chern-Simons theory with both compact and complex gauge groups.
(The extension to links is straightforward and will not be mentioned explicitly hereafter.)
This will first require a brief discussion of quantization.

%%%%%%%%%%%%%%%%%%%%%%%%%%%%%%%%%%%%%%%%%%%%%%%%%%%%%%%%%%%%%%%%%%%%%%%%%%%%%%%%%%%%%%%%%%%%%

\subsection{Quantization}
\label{sec:quant}

The basic problem of quantization begins with a pair $(\CM,\omega)$,
where $\CM$ is a manifold with symplectic structure $\omega$, called a classical ``phase space.''
Quantization takes this pair and constructs a quantum Hilbert space $\CH$.
Moreover, quantization should map the algebra of functions on $\CM$
to an algebra $\CA_{\hbar}$ of operators on $\CH$:
\be \begin{array}{ccc}
 (M,\omega) & \leadsto & \CH ~(=\text{Hilbert space}) \\
  && \raisebox{-.1cm}{\includegraphics[width=.5cm]{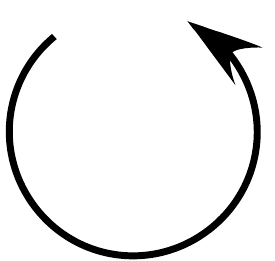}} \\
 \text{alg. of functions on $\CM$} & \leadsto & \text{alg. of operators on $\CH$} \\
 f & \mapsto & \CO_f :\CH\to\CH\,.
\end{array} \ee
The functions on $\CM$ form a Poisson algebra with respect to the usual pointwise
multiplication of functions and a Lie algebra structure $\{\bullet,\bullet\}$
induced by the symplectic structure. Quantization must map this algebra to
an associative but noncommutative algebra $\CA_{\hbar}$, such that
\be [\CO_f,\CO_g] = -i\hbar\, \CO_{\{f,g\}} + \ldots \,, \ee
where $[\bullet,\bullet]$ is the commutator of operators.
Here $\hbar$ is a parameter that is involved in the determination of $\CH$ itself as well as the algebra of operators.

Very roughly, the Hilbert space $\CH$ consists of $L^2$ sections of a complex line bundle over $\CM$ with curvature $\frac{1}{\hbar}\omega$. Locally, these sections are only allowed to depend on \emph{half} of the coordinates of $\CM$. In a standard physical setup, $\CM$ can be thought of as the space of all possible positions $x_i$ and momenta $p_i$ of particles; thus the elements of $\CH$ are functions (``wavefunctions'') that depend on either positions or moments, but not both.

In addition to the construction of $\CH$, the process of quantization must also explain how classical motions or trajectories of a physical system are associated to quantum states in $\CH$. A classical trajectory (or ``semiclassical state'') is described by a Lagrangian submanifold $\CL\subset\CM$. Being Lagrangian means that $\CL$ is middle-dimensional and $\omega|_\CL=0$. Let $\theta$ be 1-form (called a Liouville 1-form%
\footnote{There is an ambiguity in choosing $\theta$, directly related to the choice of coordinates of $\CM$ (positions versus momenta) that elements of $\CH$ are to depend on.}) %
that satisfies $\omega = d\theta$. Notice that $\theta|_\CL$ is closed. Then the Lagrangian $\CL$ is called quantizable if
\be \oint_\gamma \theta \in 2\pi \hbar\, \Z \label{thetaq} \ee
for any closed cycle $\gamma\subset\CL$. The vector (or wavefunction) $Z\in\CH$ corresponding to $\CL$ can be written as
\be Z = Z(x) = \exp\left( \frac{i}{\hbar}S_0(x) + \ldots\right)\,, \label{v} \ee
with
\be S_0(x) = \int_{x_0}^x \theta \ee
for some fixed $x_0$ and varying $x\in\CL$. Due to the quantization \eqref{thetaq}, the expression \eqref{v} is completely well-defined.

Expression \eqref{v} only defines $Z$ to leading order in $\hbar$. To find subleading corrections, it is useful to employ a complementary approach. Suppose that the Lagrangian submanifold $\CL$ is cut out by certain equations $f_i=0$ on $\CM$. Quantization promotes these functions to operators $\CO_{f_i}$ acting on $\CH$, and the vector $Z$ can also be defined as a solution to the equations
\be \CO_{f_i}\cdot Z = 0  \qquad\forall\; i\,. \ee
If the $\CO_{f_i}$ are properly quantized, then the solution to these equations will be the exact wavefunction.

\subsubsection{Methods}

The problem of quantization can be approached in many different ways.
Each approach has its advantages and disadvantages, but in the end all methods are expected to yield the same result.
The classic approach of \emph{geometric quantization} (\cf\ \cite{woodhouse-1992})
starts by defining a prequantum line bundle $L \to \CM$ with a unitary connection
of curvature $\frac{1}{\hbar}\omega$. Note that such a line bundle only exists for
\be \frac{1}{2\pi\hbar}\omega\in H^2(\CM;\Z)\,, \label{wquant} \ee
which can lead to a quantization of $\hbar^{-1}$ (\ie\ a restriction of $\hbar$ to a discrete set of values in $\C^*$).
The local choice of ``position'' versus ``momentum'' coordinates is encoded
in the choice of a set of $\frac12\dim_\R\CM$ vector fields $\CP_j$, called a polarization,
and the Hilbert space $\CH$ is then defined as the set of square-integrable,
$\CP_j$-invariant sections of $\CL$. This gives a very concrete definition of $\CH$,
although it can be very hard to show that the construction is independent of the choice of polarization.
(The problem becomes more manageable if $\CM$ is K\"ahler.)
Moreover, it is often difficult in geometric quantization
to find the full quantum expressions for operators $\CO_{f_i}$. \\

An alternative, \emph{deformation quantization} \cite{bayen-def} partially solves this latter problem.
It describes a formal $\hbar$-deformation of the ring of functions on $\CM$,
using a noncommutative product of the type
\bea
f \star_\hbar g & = & f g + \hbar \sum_{i,j} \alpha^{ij} \p_i (f) \p_j (g)
+ \frac{\hbar^2}{2} \sum_{i,j,k,l} \alpha^{ij} \alpha^{kl} \p_i \p_k (f) \p_j \p_l (g) \nonumber\\
& & + \frac{\hbar^2}{3} \left(\sum_{i,j,k,l} \alpha^{ij} \p_j (\alpha^{kl}) (\p_i \p_k (f) \p_l (g) - \p_k (f)
\p_i \p_l (g)) \right) + \ldots,
\label{starproduct}
\eea
where $\alpha = \omega^{-1}$ is the Poisson structure corresponding to the symplectic form $\omega$.
In local coordinates $\{ f,g \} = \alpha^{ij} \p_i (f) \p_j (g)$.
One important advantage of deformation quantization is that it is completely canonical
and does not require any auxiliary choices.
In particular, there is an explicit formula for the $\star_\hbar$-product \eqref{starproduct}
due to Kontsevich \cite{kont-poisson}, that allows one to express it as a sum over {\em admissible graphs},
\be
f \star_\hbar g := \sum_{n=0}^\infty \hbar^n \sum_{{{\rm graphs}\;\Gamma \atop {\rm of~order}\; n}}
w(\Gamma) B_\Gamma(f,g)\,,
\label{kproduct}
\ee
where $w(\Gamma)$ is a weight (number) assigned to $\Gamma$,
and $B_\Gamma (f,g)$ are bilinear differential operators
whose coefficients are differential polynomials, homogeneous of degree $n$
in the components of the bivector field $\alpha$ on $\CM$.
By definition, an admissible graph of order $n$ is an ordered pair of maps
$i,j : \{ 1, \ldots , n \} \to \{1, \ldots , n, L, R \}$
where neither map has fixed points and both maps are distinct at every point.
There are $n^n (n+1)^n$ such graphs.

\begin{figure}[ht]
\centering
\includegraphics[width=1.5in]{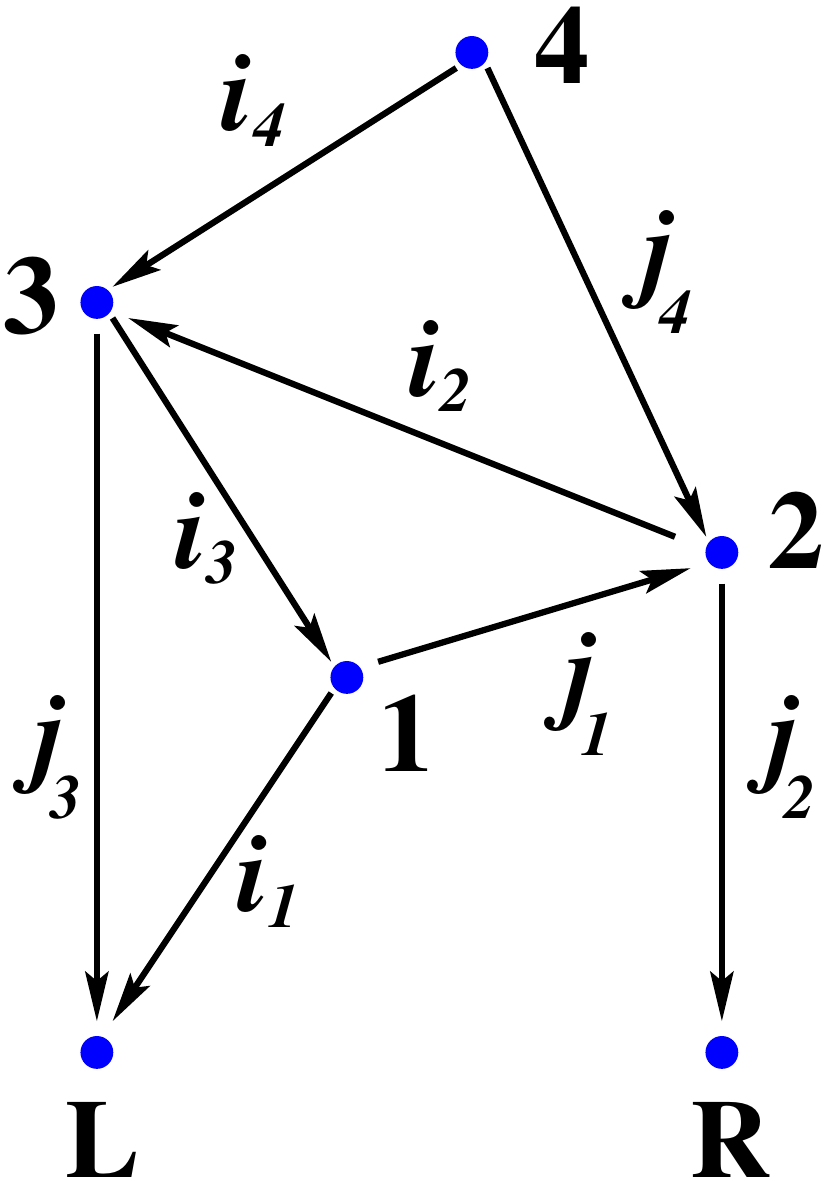}
\caption{An example of an admissible graph of order 4.}
\label{fig:Pgraph}
\end{figure}

For example, the graph of order 2 corresponding to the first term
in the second line of eq. \eqref{starproduct} has 4 vertices and 4 edges:
\be
\Gamma ~=~
\raisebox{-1.2cm}{\includegraphics[width=3.0cm]{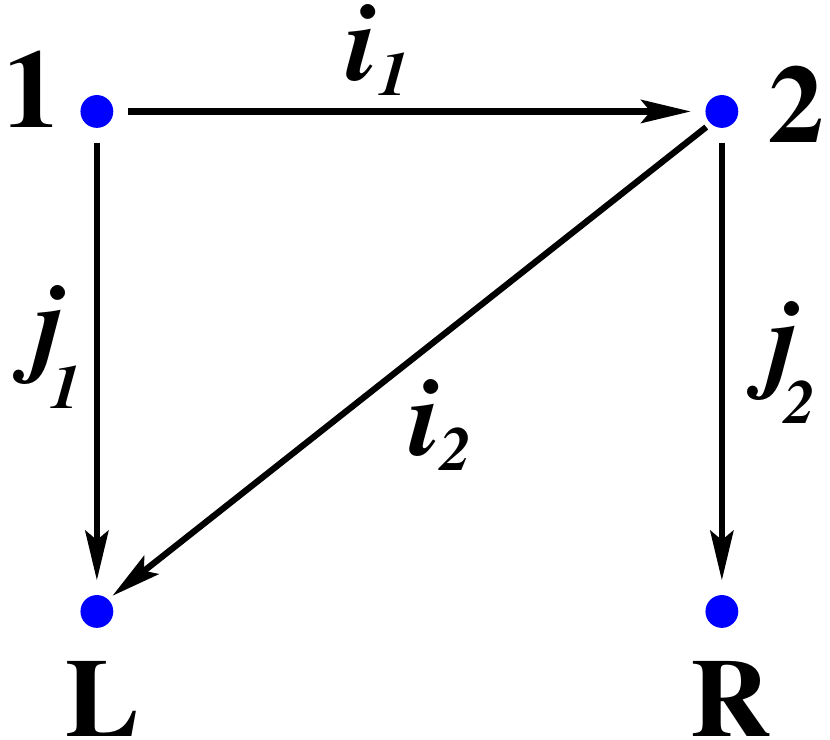}}
\qquad\qquad
\begin{array}{c}
i_1 = (1,2) \\
j_1 = (1,L) \\
i_2 = (2,L) \\
j_2 = (2,R)
\end{array} \ee
An example of a more complicated admissible graph (of order 4) is shown on Figure \ref{fig:Pgraph}.
The corresponding bidifferential operator is
\be
B_\Gamma (f,g) = \sum \alpha^{i_4 j_4} (\p_{i_3} \a^{i_1 j_1}) (\p_{j_1} \p_{j_4} \a^{i_2 j_2})
(\p_{i_2} \p_{i_4} \a^{i_3 j_3}) (\p_{i_1} \p_{j_3} f) (\p_{j_2} g)\,.
\ee

When the Poisson structure is flat, a graph with an edge ending
in a vertex other than $L$ or $R$ will have zero contribution to the sum \eqref{kproduct},
since it will involve derivatives of $\alpha$.
In this case the $\star_{\hbar}$-product \eqref{kproduct} becomes the usual Moyal product
\be
f \star_{\hbar} g (x) = \exp \left( \hbar \alpha^{ij} \frac{\p}{\p x^i} \frac{\p}{\p y^j} \right) f(x) g(y) \vert_{y=x}
\label{moyal}
\ee

Deformation quantization is a powerful method for finding the operators $\CO_{f_i}$.
It is important to stress, however, that, by itself, it does \emph{not} explain how to construct the space $\CH$
(it is not an honest quantization), and can not capture the fact that $\hbar^{-1}$ should ever be discretized. \\

A third option, \emph{brane quantization} \cite{gukov-2008},
is a marriage of geometric and deformation quantizations in a physical context.
It approaches the problem of quantization by complexifying $\CM$ and $\omega$,
and constructing a certain (secondary) topological quantum field theory on the resulting space $\CM_\C$.
It has the advantage of easily characterizing the various choices that one must make in quantization,
and provides simple geometric criteria that describe quantizable $(\CM,\omega;\hbar)$.
In this approach, the Hilbert space $\CH$ is constructed as the space
of morphisms (space of open strings),
\be
\CH = {\rm Hom} (\CB_{cc},\CB') \,,
\ee
where $\CB_{cc}$ and $\CB'$ are objects (branes)
of a certain category associated to the symplectic manifold $\CM_\C$.
Moreover, in this approach, independence of $\CH$ on various choices can be
reformulated as a problem of constructing a flat connection on the space of such choices,
which identifies the space of ground states in the secondary TQFT.
In a closely related context, this problem has been studied in the mathematical physics
literature \cite{CecottiV,Dubrovin},
and leads to a beautiful story that involves integrable systems and $tt^*$ equations.

%%%%%%%%%%%%%%%%%%%%%%%%%%%%%%%%%%%%%%%%%%%%%%%%%%%%%%%%%%%%%%%%%%%%%%%%%%%%%%%%%%%%%%

\subsubsection{Simple examples}

Let us now adapt the general statements here to some specific examples.%, building up to our main interest, Chern-Simons theory.

\subsubsection*{Harmonic oscillator}

\begin{figure}[htb]
\centering
\includegraphics[width=5in]{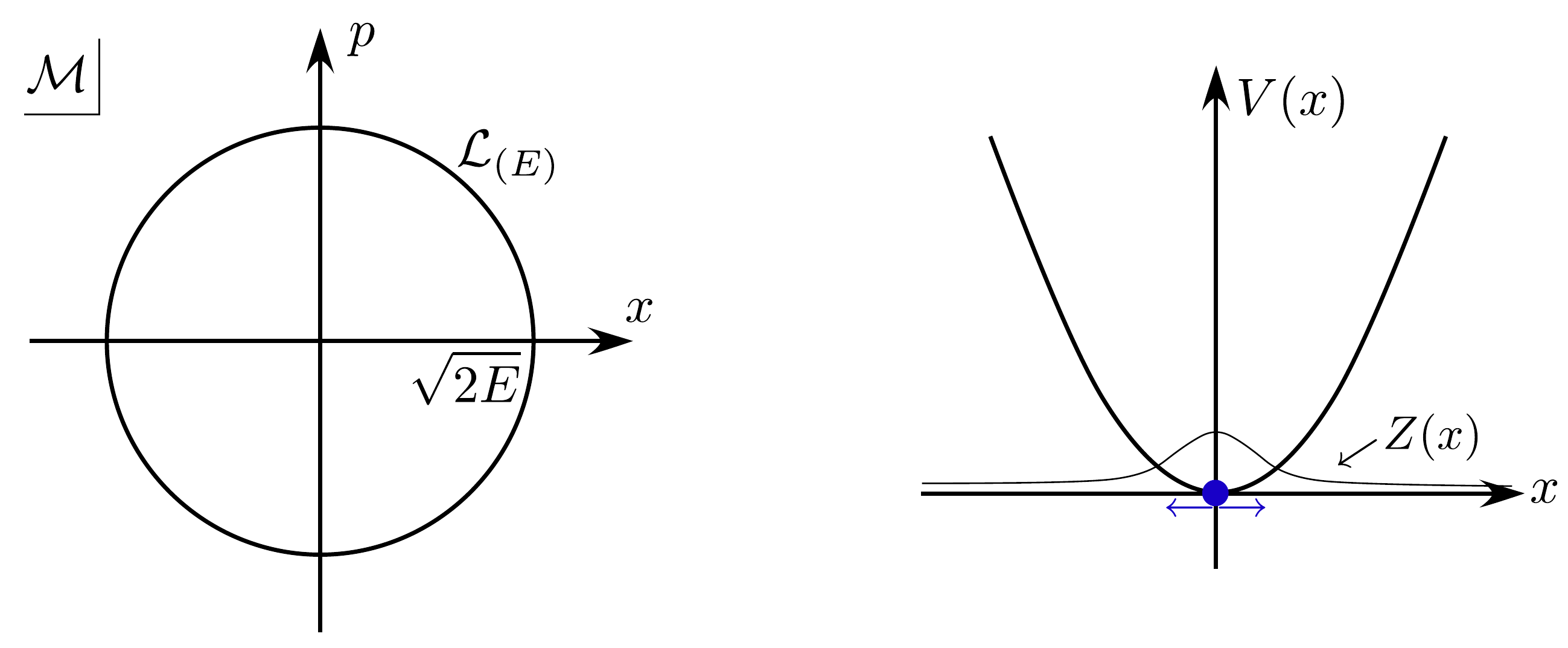}
\caption{The harmonic oscillator: potential $V(x) = \frac12 x^2$ in physical space, phase space $\CM$, a classical trajectory $\CL_{(E)}$ in phase space, and the ground state quantum wavefunction $Z(x)$.}
\label{fig:harm}
\end{figure}

The quintessential simplest nontrivial problem of quantization is the harmonic oscillator. Consider a classical system that consists of a particle moving on a line (with coordinate $x=x(t)$) with a potential energy $V = \frac12 x^2$. This is depicted in Figure \ref{fig:harm}. The total (potential + kinetic) energy of the particle at any moment of time is given by the Hamiltonian
\be H = \frac12 x^2 + \frac12 p^2\,, \ee
where classically $p = \dot x = \frac{dx}{dt}$ is the momentum. This total energy $H$ is conserved.
The classical phase space $\CM$ is just $\R^2 = \{(x,p)\}$, endowed with a symplectic structure $\omega = dp\wedge dx$. A classical trajectory with energy $H=E$ is just a circle of radius $\sqrt{2E}$ in phase space. This defines a Lagrangian submanifold $\CL_{(E)} \simeq S^1$.

Now let us quantize the system.
Since $H^2(\CM;\Z) = 0$, there is no restriction or quantization of $\hbar^{-1}$.
On the other hand, there \emph{is} a restriction on $\CL$ which quantizes the energy.
Namely, according to \eqref{thetaq},
for a Liouville 1-form $\theta$ such that $\omega = d\theta$, the integral
\be \oint_\CL \theta = \oint_{S^1} \theta = 2\pi E \ee
must be an element of $2\pi\hbar\,\Z$, implying that $E = n\hbar$ for positive $n\in\Z$.
In fact, this equation is corrected by quantum effects --- a Maslov correction in geometric quantization --- to
\be E = \hbar\Big(n+\frac12\Big)\,. \label{En} \ee
This leads to the famous result that the lowest possible energy of a quantum harmonic oscillator (at $n=0$) is nonzero.

Suppose we choose a polarization $\pd/\pd p=0$, and a corresponding Liouville 1-form $\theta = p\,dx$. The Hilbert space $\CH$ can simply be identified as $L^2(\R) \sim \{\text{functions of $x$}\}$, on which the functions $x$ and $p$ act as operators
\be \hat{x}:=\CO_x = x\,, \qquad \hat{p}:=\CO_p = -i\hbar\frac{d}{dx}\,. \ee
In this case, the exact quantum expression for the Hamiltonian is
\be \CO_H = \frac12(\hat{x}^2+\hat{p}^2)\,. \ee
It is then easy to find the quantum wavefunctions corresponding to classical states $\CL_{(E)}$. From \eqref{v}, we find a leading contribution
\be Z(x) \simeq \exp\left(\frac{i}{\hbar}\int_0^x \theta\right) = \exp\left(\frac{i}{\hbar}\int_0^x\sqrt{2E-x^2}\,dx\right)
 \simeq \exp\left(-\frac{1}{2\hbar}x^2+ \ldots\right)\,.
\ee
Since the Lagrangian $\CL_{(E)}$ is defined classically by $H-E=0$,
the complete expression for $Z(x)$ can be obtained by solving the operator equation $(\CO_H-E)Z=0$.
This eigenvalue equation has square-integrable solutions only for the quantized energies \eqref{En};
for example, at the ground state energy $E = \hbar/2$,
the exact solution is $Z(x) = \exp\left(-\frac{1}{2\hbar}x^2\right)$.

%%%%%%%%%%%%%%%%%%%%%%%%%%%%%%%%%%%%%%%%%%%%%%%%%%%%%%%%%%%%%%%%%%%%%%%%%%%%%%%%%%%%%%%%%%%%%%%%%%%%%

\subsubsection*{Representations of Lie groups}

Another famous application of quantization
is the construction of unitary representations of Lie groups by quantization of coadjoint orbits.
A basic premise of this approach, also known as the {\em orbit method},
is that coadjoint orbits come equipped with a natural symplectic structure
(the Kostant-Kirillov-Souriau symplectic structure), therefore
providing interesting examples for quantization.

Continuing with our default notations in these notes,
we use $G$ for a compact Lie group (that we usually assume to be simple),
$G_{\C}$ for its complexification,
and $G_{\R}$ for some real form of the complex group $G_{\C}$ (that may be equal to $G$).
We denote by $\fg_{\R}$ the Lie algebra of $G_{\R}$, and similarly for $G$ and $G_{\C}$.
Given an element $\lambda \in \fg_{\R}^*$ (the highest weight of the desired unitary representation $R_\lambda$)
one constructs $\CM = G_{\R} \cdot \lambda$
as the coadjoint orbit of $G_{\R}$ in $\fg_{\R}^*$ passing through $\lambda$.

In the case of compact groups, the phase space $\CM$ is compact
and its quantization leads to a finite-dimensional Hilbert space $\CH$
as the space of the unitary representation $R_\lambda$.
This is the statement of the Borel-Bott-Weil theorem.
Moreover, the condition $\frac{1}{2\pi\hbar}\omega\in H^2(\CM;\Z)$ that ensures the existence of a prequantum line bundle
becomes equivalent to the condition that $\lambda$ be an element of the weight lattice $\Lambda_w\subset\fg^*$.

\begin{comment}
,
so that one finds
\be
\label{bwbspace}
\CH = H^0_{\bar \p} (\CM, \CO (\lambda + \rho)) \,.
\ee
This is precisely the statement of the Borel-Weil-Bott theorem
which asserts that, for each $\lambda \in H^2(\CM;\Z)$,
an irreducible representation of $G_{\R} = G$
with highest weight $\lambda$ can be realized as the space
of holomorphic sections of $\CO (\lambda + \rho)$, as in \eqref{bwbspace}.
\end{comment}

As a very simple illustration, consider the group $SU(2)$.
In this case, a non-trivial coadjoint orbit is topologically equivalent to the flag manifold
\be SU(2)/U(1) \simeq \P^1 \,. \ee
Letting $\omega$ be the unit volume form on $\P^1$, we see that $(\CM,\omega)$ is quantizable for
\be \hbar^{-1} = 2\pi\lambda\,, \qquad \lambda\in\Z_{(+)}\,. \ee
The prequantum line bundle with curvature $\hbar^{-1}\omega$ is simply $\CO(\lambda)\to \P^1$.
Choosing a holomorphic polarization, so that $\CH$ is defined as the space of holomorphic sections
of $\CO(\lambda)$, we see that $\dim\CH = \lambda+1$.
The Hilbert space is precisely the space of the $(\lambda+1)$-dimensional representation of $SU(2)$.

Similarly, some infinite-dimensional representations,
such as unitary principal series representations of $SL(n,\C)$ or $SL(n,\R)$, can be described as quantized orbits.
Nevertheless, there remain some outstanding puzzles:
there exist unitary representations that don't appear to correspond to orbits, and, conversely, there are real orbits that don't seem to correspond to unitary representations.
An example of first kind occurs even in the basic case
of the real group $G_{\R} = SL(2,\R)$ and the complementary series representations.
To illustrate the second phenomenon, one can take $G_{\R}$ to be
a real group of Cartan type $B_N$, {\it i.e.} $G_{\R} = SO(p,q)$ with $p+q=2N+1$.
The minimal orbit $\CO_{{\rm min}}$ of $B_N$ is a nice symplectic manifold
of (real) dimension $4N-4$, for any values of $p$ and $q$.
On the other hand, the minimal representation of $SO(p,q)$
exists only if $p \le 3$ or $q \le 3$ \cite{vogan-singular-unitary}.
Both of these issues can be resolved in the brane quantization approach \cite{gukov-2008},
at the cost of replacing classical geometric objects (namely, coadjoint orbits)
with their quantum or ``stringy'' analogs (branes).
In particular, in the case of $B_N$ one finds that,
while the minimal orbit exists for any values of $p$ and $q$,
the corresponding brane exists only if $p \le 3$ or $q \le 3$.
(In general, the condition is that the second Stieffel-Whitney class
$w_2(\CM) \in H^2(\CM;\Z_2)$ must be
a mod 2 reduction of a torsion class in the integral cohomology of $\CM$.)

%%%%%%%%%%%%%%%%%%%%%%%%%%%%%%%%%%%%%%%%%%%%%%%%%%%%%%%%%%%%%%%%%%%%%%%%%%%%%%%%%%%%%%%%%%%%%%

\subsection{Chern-Simons theory}

Finally, we arrive at our goal, Chern-Simons theory. Let us recall for a second why we began discussing quantization in the first place. In Section \ref{sec:cut}, we reviewed how partition functions in TQFT could be obtained by cutting and gluing three-manifolds. We explained that the partition function of a manifold with a knot is equivalent to the partition function of the corresponding knot complement, projected onto appropriate boundary conditions in $\CH(T^2)$ as in \eqref{uR}. To make complete sense of this, however, and to actually calculate partition functions, we must understand what $\CH(T^2)$ really is. Using Section \ref{sec:quant} we are finally in a position to do so.

\subsubsection{Quantization of Chern-Simons theory}
\label{sec:TQFT_CS}

Consider Chern-Simons theory with gauge group $G$ --- either compact or complex --- on a knot complement $M = \widetilde{M}\bs K$, with $\pd M = T^2$. The phase space $\CM$ associated to $T^2$ is simply the space of flat $G$-connections on $T^2$, modulo gauge equivalence. Since a flat connection is completely determined by the conjugacy classes of its holonomies, we have
\begin{align} \CM &= \CM_{{\rm flat}} (G;T^2) \\
 &= \left\{ \text{representations :}\;\;\pi_1(T^2)\to G \right\}/\text{conjugation}\,.
\end{align}
The fundamental group $\pi_1(T^2)\simeq \Z\oplus\Z$ is abelian, generated by the meridian and longitude of the torus. The holonomies along these loops can therefore be simultaneously diagonalized%
\footnote{If $G$ is not compact, there may be elements that are not so diagonalizable, but they form lower-dimensional components of $\CM$ which should not be considered in the quantization.} %
into the maximal torus $\mb{T}\subset G$. Coordinates on $\CM$ are then given by the $2r$ independent eigenvalues $(m_1,...,m_r)$ and $(\ell_1,...,\ell_r)$ of the meridian and longitude holonomies, where $r$ is the rank of $G$. We must also divide by the Weyl group $\CW$ of $G$, which simultaneously permutes both sets of eigenvalues, to obtain
\be \CM \simeq (\mb{T}^r\times \mb{T}^r)/\CW = \mb{T}^{2r}/\CW\,. \ee
For example, for a compact group $G = SU(n)$
the phase space is $\CM = (S^1)^{2(n-1)}/S_n$,
where $S_n$ is the symmetric group on $n$ elements.
Similarly, for its complexification $G_\C = SL(n,\C)$, the phase space is $\CM = (\C^*)^{2(n-1)}/S_n$.
In general, ignoring subtleties in high codimension that are not pertinent to quantization,
the relation between compact and complex phase spaces can be described as
\be \CM_{{\rm flat}} (G_\C;\Sigma) = \left[\CM_{{\rm flat}} (G;\Sigma)\right]_\C \simeq T^*\CM_{{\rm flat}} (G;\Sigma)\,. \ee
(In particular, the last relation is only a birational equivalence.)

%%%%%%%%%%%%%%%%%%%%%%%%%%%%%%%%%%%%%%%%%%%%%%%%%%%%%%%%%%%%%%%%%%%%%%%%%%%%%%%%%%%%%

\subsubsection*{Compact theory}

In order to quantize $\CM$, we need a symplectic structure.
In \emph{compact} Chern-Simons theory, it is given by
\be \omega = \frac14 \int_{T^2} \Tr\big[\delta\CA\wedge \delta\CA\big]\,. \label{w0} \ee
This can be expressed more concretely in coordinates $\{m_i,\ell_i\}=\{e^{u_i},e^{v_i}\}$ as
\be \omega = \sum_i d\log m_i\wedge d\log \ell_i = \sum_idu_i\wedge dv_i\,. \ee
The holonomy variables $u_i$ and $v_i$ function as ``positions'' and ``momenta,'' respectively.
Now, the parameter $\hbar=i\pi/k$ that appeared naturally in the discussion of the volume conjecture in Section \ref{sec:gen} is rescaled from the standard geometric quantization parameter $\hbar$ of Section \ref{sec:quant} by a factor of $i$. In terms of $k$, the quantization condition \eqref{wquant} simply takes the form $k\in\Z$. The integer $k$ is identified as the Chern-Simons level, modulo the shift mentioned in Footnote \ref{foot:k}.

The last ingredient we need to describe the Hilbert space $\CH$ is a choice of polarization. For clarity, let us take $G=SU(2)$ to be of rank one, and let us choose the polarization $\pd/\pd v = 0$, so that $\CH_{T^2}$ essentially consists of periodic and Weyl-invariant functions of $u$, $f(u) = f(u+2\pi i) = f(-u)$. Being somewhat more careful, and thinking of these not as functions but as sections of the line bundle with curvature $\frac{k}{\pi}\omega$, one finds that the simultaneous periodicity in the momentum $v$ and the position $u$ restricts $u$ to take values in $\frac{i\pi}{k}\Z$. Therefore, a function $f(u)$ only takes nonzero values at $k+1$ distinct points $u=0,\frac{i\pi}{k},\frac{2 i\pi}{k}...,i\pi$, and the space $\CH_{T^2}$ is finite-dimensional. For general compact semi-simple $G$, the Hilbert space $\CH$ takes the form \cite{EMSS, axelrod-witten}
\be \CH_{T^2} \simeq \frac{\Lambda_w}{\CW\ltimes k\Lambda_r}\,. \ee
where $\Lambda_w,\,\Lambda_r$ are the weight and root lattices of $G$. In other words, $\CH_{T^2}$ is the set of weights (hence representations) in a level-$k$ affine Weyl chamber.

Given a ``wavefunction'' $Z(M;u;\hbar)\in\CH_{T^2}$ associated to the knot complement $M=\widetilde{M}\bs K$
(with $\pd M=T^2$), the partition function $Z^G(\widetilde{M},K;R_\lambda;\hbar)$
for $K\in \widetilde{M}$ colored by representation $R_\lambda$ is simply given by
evaluating $Z(M;u;\hbar)$ at $u=i\pi\frac{\lambda^*+\rho^*}{k}$ as in \eqref{uR}.
For example, in the case of $SU(2)$ theory, we evaluate $Z(M;u;\hbar)$ at $u = i\pi N/k$
(and normalize by the partition function of $S^3$) to find the colored Jones polynomial $J_N(K;q)$.
The single wavefunction $Z(M;u;\hbar)$ in $\CH_{T^2}$
comprises the entire family of colored Jones polynomials $J_N(K;q)$, $N\in\Z$.

How is such a wavefunction obtained in the compact theory? For any three-manifold $M$, there is a Lagrangian submanifold $\CL \subset \CM$ corresponding to the semi-classical ``state'' $M$. This manifold $\CL$ is simply defined as the set of flat connections on $T^2$ that can extend to a flat connection on all of $M$. It is the so-called $G$-character variety of $M$ and can be described by a set of polynomial equations in the eigenvalues $\ell_i$ and $m_i$:
\be \CL\,:\quad A_j(\ell,m) = 0\,. \ee
Depending on whether we restrict to $\ell,m\in S^1$ or $\ell,m\in\C^*$, these same equations describe flat $G$ or $G_\C$ connections. In the rank-one case, there is just a single equation, the A-polynomial of the knot complement. Upon quantization, the functions $A_j$ get mapped to quantum operators
\be \widehat{A}_j(\widehat{\ell},\widehat{m},q=e^{\frac{2\pi i}{\hbar}}) := \CO_{A_j}\,,\ee
where $\widehat{\ell}_i :=\CO_{\ell_i} = e^{\widehat{v}_i}$
and $\widehat{m}_i :=\CO_{m_i} =e^{\widehat{u}_i}$ act on $\CH_{T^2}$ as
\be \widehat{\ell}_i\,Z(u) = Z(u_i+\hbar)\quad\text{(shifting only $u_i$)}\,,\qquad \widehat{m}_i\,Z(u) = e^{u_i} Z(u)\,. \ee
In terms of the colored Jones polynomial $J_N(K,q)$,
this means $\widehat{\ell}J_N(K,q) = J_{N+1}(K,q)$ and $\widehat{m}J_N(K,q) = q^{N/2}J_N(K,q)$.
The wavefunction $Z(M;u)$ must satisfy \cite{gukov-2003, DGLZ}
\be \widehat{A}_j\,Z(M;u) = 0 \qquad \forall\,j\,, \label{Acpt} \ee
which leads to a set of recursion relations on polynomial invariants of the knot $K$.
In the mathematical literature, such a recursion relation for the colored Jones polynomial
(\ie\ in the case of $G=SU(2)$) is known as the {\em AJ conjecture} \cite{garoufalidis-2004, Gar-Le}
(also \cf\ \cite{Gar-twist}).

%%%%%%%%%%%%%%%%%%%%%%%%%%%%%%%%%%%%%%%%%%%%%%%%%%%%%%%%%%%%%%%%%%%%%%%%%%%%%%%%%%%%%%%%

\subsubsection*{Complex theory}

Now, we would like to relate partition functions in Chern-Simons theory with compact gauge group $G$ to Chern-Simons theory with complex gauge group $G_\C$. In the case of complex gauge group, the phase space is $\CM=\CM_{{\rm flat}} (G_\C;T^2)=\big((\C^*)^r\times(\C^*)^r\big)/\CW$, and the full symplectic structure induced by Chern-Simons theory is
\be \omega = \frac{\tau}{2}\,\omega_0 + \frac{\tilde{\tau}}{2}\,\ol{\omega_0}\,, \label{wcx} \ee
with $\omega_0 = \frac14\int_{T^2} \Tr(\delta \CA\wedge\delta\CA)$ as in \eqref{w0}. The connection $\CA$ now takes values in $\fg_\C$, and a priori there are two independent coupling constants $\tau$ and $\tilde{\tau}$. These are the analog of the level $k$ in the compact theory; we include them here in the definition of $\omega$. Since $\CM$ is noncompact, the quantization condition \eqref{wquant} is less restrictive, only fixing $\tau+\tilde{\tau}\in \Z$.

The noncompactness of $\CM$ changes the nature of the Hilbert space $\CH$ --- as in the case of the harmonic oscillator, it is no longer finite-dimensional. Choosing a polarization $\pd/\pd v=0$, we can effectively take $\CH$ to  consist of Weyl-invariant square-integrable functions $f(u,\bar{u}) \in L^2((\C^*)^{r})$. However, the fact that \eqref{wcx} is a simple sum of holomorphic and antiholomorphic pieces means that at a perturbative level any wavefunction $Z(M;u)\in\CH$ will \emph{factorize} into holomorphic and antiholomorphic components. Put more concretely, the exact wavefunction $Z(M;u,\bar{u};\hbar=\frac{2\pi i}{\tau},\tilde{\hbar}=\frac{2\pi i}{\tilde{\tau}})$ corresponding to complex Chern-Simons theory on the knot complement $M$ can be written as \cite{DGLZ, Wit-anal}
\be Z(M;u,\bar{u};\hbar,\tilde{\hbar}) = \sum_{\alpha,\bar{\alpha}} n_{\alpha,\bar{\alpha}}\, Z^\alpha_{\rm pert}(M;u;\hbar)\, \ol{Z}^{\bar{\alpha}}_{\rm pert}(M;\bar{u};\tilde{\hbar})\,, \label{sumGc}\ee
for some coefficients $n_{\alpha,\bar{\alpha}}$, where, as $\hbar\to 0$, each component $Z^\alpha_{\rm pert}(M;u)$ can be expressed as a perturbative series
\be  Z^\alpha_{\rm pert}(M;u) = \exp\left( -\frac{1}{\hbar}S_0(u) +\frac{\delta}{2}\log\hbar+S_1(u)+ \hbar\,S_1(u)+\ldots\right)\,,\qquad \hbar = \frac{2\pi i}{\tau}\,\,. \label{Zcxpert} \ee

Each partition function $Z^\alpha_{\rm pert}(M;u;\hbar)$ corresponds to complex Chern-Simons theory in the background of a fixed flat connection on $M$ that has meridian holonomy eigenvalues $m=e^u$. The set of such connections, labelled by $\alpha$, is nothing but the (finite) set of solutions
$\{v^\alpha(u)\}$ (mod $2\pi i$) to the equations
\be A_j(\ell,m) = 0 \ee
at fixed $m=e^u$. In the case of $SL(2,\C)$ theory, one of these flat connections is the geometric one, corresponding to a hyperbolic metric on $M$.

Since the complex phase space $\CM$ is just the complexification of the phase space of the compact theory, the quantization of the functions $A_j(\ell,m)$ is formally identical to the quantization in the compact case. In other words, the operators $\CO_{A_j}=\widehat{A}_j(\hat{\ell},\hat{m},q=e^{2\hbar})$ are identical to those of the compact theory. Every component $Z^\alpha_{\rm pert}(M;u;\hbar)$ must therefore satisfy \cite{gukov-2003, DGLZ}
\be \widehat{A}_j(\widehat{\ell},\widehat{m},e^{2\hbar})\,Z^\alpha_{\rm pert}(M;u;\hbar) = 0\,\qquad \forall\,j,\alpha\,,  \label{Acx} \ee
with $\widehat{\ell}_iZ(M;u;\hbar) = Z(M;u_i+\hbar;\hbar)$ (in other words $\widehat{v}_i=\hbar\pd_{u_i}$) and $\widehat{m}_iZ(M;u;\hbar) = u^{u_i}Z(M;u;\hbar)$.
In particular, at leading order in $\hbar$, we can write
\be Z^\alpha_{\rm pert}(M;u;\hbar) = \exp\left(-\frac{1}{\hbar}\int_{\gamma_\alpha} \theta +\ldots\right)\,, \label{Zcts} \ee
where $\theta \sim - \sum_iv_i\,du_i$ is a Liouville 1-form and $\gamma_\alpha$
is a path on the complex variety $\CL = \{ A_j=0\}$
ending at the point $(e^{v^\alpha(u)},e^u)$, as in Figure \ref{fig:lm}b.
Now that $u$ is a continuous parameter in the complex theory, this integral expression makes complete sense.

%%%%%%%%%%%%%%%%%%%%%%%%%%%%%%%%%%%%%%%%%%%%%%%%%%%%%%%%%%%%%%%%%%%%%%%%%%%%%%%%%%%%%%%%%%%%%%%%%%%%%%%%

\subsubsection{Synthesis}
\label{sec:path}

It is fairly clear from the above discussion of quantization that there should be a relation between the partition function for Chern-Simons theory with compact gauge group $G$ and the partition function for Chern-Simons theory with complex gauge group $G_\C$. Essentially the same equations \eqref{Acpt} and \eqref{Acx} define the two partition functions --- though in one case they are difference equations and in the other they are differential equations. This relation was developed in \cite{gukov-2003, DGLZ}, and was recently explained very concretely in \cite{Wit-anal} in terms of analytic continuation.

Algebraically, there may be several solutions to the difference equations \eqref{Acpt} of the compact theory.
Let us label them as $Z^{\alpha}_G(M;u)$.
The exact partition function of the compact theory
(\ie\ the colored Jones polynomial for $G=SU(2)$)
is given as a linear combination
\be Z^G(M;u;\hbar) = \sum_\alpha n_\alpha Z^\alpha_G(M;u;\hbar)\,. \label{sumG} \ee
The $\hbar\to 0$ asymptotics of each component in this sum
are then governed by the corresponding solution $Z^\alpha_{\rm pert}(M;u)$
to the differential equation \eqref{Acx}, written in the form \eqref{Zcxpert}.
These are holomorphic pieces of the $G_\C$ partition function.
The physical statement of the volume conjecture for $SU(2)$
is that the component of the sum \eqref{sumG} with the \emph{dominant}
leading asymptotics corresponds to the $SL(2,\C)$
partition function $Z^{\alpha=\rm hyp}_{\rm pert}(M;u)$
around the hyperbolic $SL(2,\C)$ flat connection.
Of all the flat $SL(2,\C)$ connections, this has the largest volume
in a neigborhood of the complete hyperbolic point $u=i\pi$.
Therefore, \emph{if} the solution $Z^{\rm hyp}_G(M;u)$
of the difference equations contributes to the colored Jones polynomial in \eqref{sumG},
it will have the dominant asymptotic. One must simply assure that
\be \text{Physical volume conjecture\,:}\qquad  n_{\rm hyp}\neq 0\,.\hspace{1in} \ee

For higher-rank groups, it is again clear that the overall asymptotics of $Z^G(M;u)$ will be controlled by the flat $G_\C$ connection with the largest volume that makes a corresponding contribution to \eqref{sumG}. One may expect by comparison to $SU(2)$ theory that the connection with the largest overall volume (the analog of the hyperbolic flat connection) in fact contributes and dominates. This has yet to be explored. \\

The expansions \eqref{sumG} and \eqref{sumGc} for compact and complex Chern-Simons theory, and the relation between them, were explained in \cite{Wit-anal} using analytic continuation of the Chern-Simons path integral. The path integral provides yet another method for quantizing a topological quantum field theory, with its own inherent advantages. Let us finish by saying a few words about this.

The path integral for compact Chern-Simons theory takes the form
\begin{align} Z^G(M;u;\hbar)  &= \int \CD\CA_{(u)} \exp\left(\frac{ik}{4\pi}I_{CS}(\CA)\right) \label{Pcpt} \\
&= \int \CD\CA_{(u)} \exp\left(-\frac{1}{4\hbar}I_{CS}(\CA)\right)\,, \nno
\end{align}
where $I_{CS}(\CA) = \int_M \Tr\big(\CA d\CA+\frac23\CA^3\big)$ is the Chern-Simons action as in \eqref{ICS} and $k\in\Z$.
The integral is over all $G$-connections on $M=\widetilde{M}\bs K$, modulo gauge equivalence, with fixed holonomy eigenvalues $e^u$ at the meridian of $K$. (In order to obtain a nonzero answer, $u/\hbar\sim \lambda^*+\rho^*$ must be integral.) For $\fg$-values connections $\CA$, the action $I_{CS}(\CA)$ is real. Therefore, for $k\in\Z$, the integral \eqref{Pcpt} is oscillatory and can be calculated by appropriately regulating the oscillations as $\CA\to\infty$. In \cite{Wit-anal}, however, the problem was posed of analytically continuing to $k\in \C$. Roughly speaking, to accomplish this one must also complexify the gauge connection $\CA$ so that it is $\fg_\C$-valued. When $k\in\R$, the integral \eqref{Pcpt} is then interpreted as a holomorphic contour integral along the real subspace in the space of complex connections. As $k$ is pushed away from the real line, this integration contour must also move. In general, the appropriate integration contour for $k\in\C$ is a sum of contours going through the various saddle points of the complexified action $I_{CS}(\CA)$. Each saddle point is a flat $G_\C$ connection, and an expression of the form \eqref{sumG} results.

For complex Chern-Simons theory, the procedure is quite similar. The path integral is
\begin{align} Z^{G_\C}(M;u;\hbar,\tilde{\hbar}) &=
 \int \CD\CA_{(u)}\CD\ol{\CA}_{(\bar{u})}
   \exp\left( \frac{i\tau}{8\pi} I_{CS}(\CA)+ \frac{i\tilde{\tau}}{8\pi} I_{CS}(\ol{\CA}) \right)  \label{Pcx} \\
 &= \int \CD\CA_{(u)}\CD\ol{\CA}_{(\bar{u})}
   \exp\left( -\frac{1}{4\hbar} I_{CS}(\CA)- \frac{1}{4\tilde{\hbar}} I_{CS}(\ol{\CA}) \right)\,, \nno
\end{align}
for a $\fg_\C$-valued $G_\C$-connection $\CA$.
The integrand is initially well-defined when $\tau+\tilde{\tau}\in\Z$, and leads to a convergent oscillatory integral when the exponent is imaginary --- \ie\ for $\tilde{\tau}=\bar{\tau}$. In order to analytically continue to independent $\tau,\tilde{\tau}\in\C$, one must treat $\CA$ and $\ol{\CA}$ as independent connections and again complexify \emph{each} of them.  One then deforms the contour of integration away from the ``real'' subspace when $\tilde{\tau}\neq\bar{\tau}$, and writes the resulting contour as a sum over \emph{pairs} of saddle points for $\CA$ and $\ol{\CA}$. Since $(\fg_\C)_\C \simeq \fg_\C\times \fg_\C$, however, these are just pairs of saddle points of flat $G_\C$-connections. An expression of the form \eqref{sumGc} results:
\be
\label{sumGG}
Z^{G_\C}(M;u;\hbar,\tilde{\hbar}) = \sum_{\alpha,\bar{\alpha}}n_{\alpha,\bar{\alpha}} Z^{\alpha}_{G_\C}(M;u;\hbar) \ol{Z^{\alpha}_{G_\C}}(M;\bar{u};\tilde{\hbar})\,.
\ee
The functions $Z^{\alpha}_{G_\C}(M;u;\hbar)$ and $Z^{\alpha}_{G}(M;u;\hbar)$ here and in \eqref{sumG} should be \emph{identical}, since they both correspond to $G_\C$ connections.

In \cite{Wit-anal}, it is explained how the coefficients $n_\alpha$
and $n_{\alpha,\bar{\alpha}}$ may be calculated for specific examples,
like the trefoil and figure-eight knot complements.
As expected, the coefficient of the hyperbolic component ``$\alpha={\rm hyp}$''
of the $SU(2)$ partition function is nonzero, leading to another demonstration of the volume conjecture.

The careful reader may still be wondering why it is only the growth of the colored Jones polynomial
at \emph{non}rational $N/k$ that shows exponential behavior.
The answer comes from a final subtlety in the analytic continuation
of the path integral: for $k\notin\Z$, the sum \eqref{sumG}
can have \emph{multiple} contributions from the same flat connection,
differing by a multiplicative factor $e^{2\pi i k}$. (If analytically continuing in $N$ as well,
factors of $e^{\frac{2\pi i u}{\hbar}}$ may also arise.)
This behavior originates from the fact that $\exp\big(\frac{ik}{4\pi}I_{CS})$
is not completely gauge-invariant when $k\notin\Z$.
For example, in the case of the figure-eight knot,
the actual hyperbolic contribution to \eqref{sumG} goes like
\be \big(e^{i\pi k}-e^{-i\pi k})Z^{\rm hyp}_{G_\C}(M;u;\hbar)\,, \ee
which vanishes at $k\in\Z$, leading to polynomial rather than
exponential growth of $J_N(K;q)$ for $u/i\pi \sim N/k\in \Q$.
It is expected that this feature is fairly generic for hyperbolic knots.

\vspace{1cm}
\noindent\textbf{Acknowledgements} \\
We would like to thank Edward Witten, Don Zagier, and Jonatan Lenells
for enlightening discussions on subjects considered in these notes.
We would also like to thank the organizers of the workshops
{\it Interactions Between Hyperbolic Geometry, Quantum Topology, and Knot Theory}; {\it Chern-Simons Gauge Theory: 20 years after}; and {\it Low Dimensional Topology and Number Theory II}; and Columbia University, the Hausdorff Center for Mathematics, and the University of Tokyo, respectively, for their generous support,
accommodations, and collaborative working environment.
The work of SG is supported in part by DOE Grant DE-FG03-92-ER40701,
in part by NSF Grant PHY-0757647, and in part by the Alfred P. Sloan Foundation.
Opinions and conclusions expressed here are those of the authors
and do not necessarily reflect the views of funding agencies.

\bibliographystyle{JHEP_TD}
\bibliography{VCReview}

\end{document}